\documentclass[11pt,oneside,reqno]{amsart}

\usepackage[hmargin=1.in,vmargin=1.0in]{geometry}
\usepackage{enumerate}
\usepackage{bm}
\usepackage{amssymb}
\usepackage{outlines}
\usepackage{epsfig}  		
\usepackage{epic,eepic}       
\usepackage{dutchcal} 
\usepackage{enumerate} 
\usepackage{bbm}
\usepackage{amsthm}
\usepackage[latin9]{inputenc}
\usepackage{babel}
\usepackage{booktabs}
\usepackage{subfigure}
\usepackage[colorlinks=true,citecolor=blue,urlcolor=blue]{hyperref}
\usepackage{amsmath}
\usepackage{amssymb,amsfonts}
\usepackage{amsfonts}
\usepackage{blindtext}
\usepackage{tabularray}
\UseTblrLibrary{amsmath, booktabs}

\usepackage{pdflscape}
\usepackage{rotating}
\usepackage{adjustbox}

\newtheorem{thm}{Theorem}

\newcommand{\At}{\tilde{\mathbf{A}}}
\newcommand{\bt}{\tilde{\mathbf{b}}}
\newcommand{\Ut}{\tilde{\mathbf{U}}}
\newcommand{\yt}{\tilde{y}}

\newcommand{\UU}{\mathbf{U}}
\newcommand{\RR}{\mathbb{R}}
\newcommand{\CC}{\mathbb{C}}
\newcommand{\AAA}{\mathbf{A}}
\newcommand{\RRR}{\mathbf{R}}
\newcommand{\ee}{\mathbf{e}}
\newcommand{\id}{\mathbf{I}}
\newcommand{\LR}{\mathbf{L}_R}
\newcommand{\LM}{\mathbf{L}_M}
\newcommand{\MM}{\mathbf{M}}

\newcommand{\diag}{\mathop{\mathrm{diag}}}
\newcommand{\OC}{\mathcal{O}}
\newcommand{\bb}{\mathbf{b}}
\newcommand{\cc}{\mathbf{c}}
\newcommand{\RE}{\mathrm{Re}}

\usepackage[T1]{fontenc}
\usepackage{array,booktabs}
\usepackage{siunitx}
 \usepackage{etoolbox}
\AtBeginEnvironment{tabular}{\small}

\begin{document}
\title{Very High-Order A-stable Stiffly Accurate Diagonally Implicit Runge-Kutta Methods with Error  Estimators}
\author{Yousef Alamri \and David I. Ketcheson}
\thanks{Computer, Electrical, and Mathematical Sciences \& Engineering Division, King Abdullah University of Science and Technology, Thuwal 23955, Saudi Arabia}
\email{yousef.alamri@kaust.edu.sa, 
 david.ketcheson@kaust.edu.sa} 
\date{\today}

\begin{abstract}
A numerical search approach is used to design high-order diagonally implicit Runge-Kutta (DIRK) schemes equipped with embedded error estimators, some of which have identical diagonal elements (SDIRK) and explicit first stage (ESDIRK).
In each of these classes, we present new A-stable schemes of order six (the highest order of
previously known A-stable DIRK-type schemes) up to order eight. For each order,
we include one scheme that is only A-stable as well as schemes that are L-stable, stiffly accurate, and/or have stage order two. The latter types require more stages, but give better convergence rates for differential-algebraic equations (DAEs), and those which have stage order two give better accuracy for moderately stiff problems. The development of the eighth-order schemes requires, in addition to imposing A-stability, finding highly accurate numerical solutions for
a system of 200 equations in over 100 variables, which is accomplished via a combination of global and local optimization strategies.
The accuracy, stability, and adaptive stepsize control of the schemes are demonstrated on diverse problems.
\end{abstract}

\maketitle

\section{Introduction}
This work primarily concerns the numerical approximation of several types of initial value problems (IVPs). One such task is the approximation of the solution $y:[t_0,T] \to \RR^m$ to the system of ordinary differential equations (ODEs) 
    \begin{equation} \label{IVP}
        y'(t) = f(t,y), \ \ \ y(t_0) = y_0,
    \end{equation}
for sufficiently smooth function $f:[t_0,T] \times \RR^m \to \RR^m$. Systems of the form $\eqref{IVP}$ may also arise from discretizations of partial differential equations (PDEs), for instance, if the method of lines approach is used for spatial discretization. One class of numerical methods for approximating the solution to $\eqref{IVP}$ is the Runge-Kutta (RK) family which belongs to the one-step multi-stage class of numerical methods. An $s$-stage RK method generates a sequence of approximations $y_n \approx y(t_n)$ at discrete points $t_1,\dots,t_n$ as well as a local error estimate $y_n - \widehat{y}_n$  as follows:
\begin{subequations} \label{RK}
\begin{align}
    u_i &= y_n + \Delta t \sum_{j=1}^s a_{i,j} f(t_n + c_j\Delta t, u_j), \ \ \ \ i = 1,\dots,s \label{RK1}\\
    y_{n+1} &= y_n + \Delta t \sum_{i=1}^s b_i f(t_n+c_i\Delta t,u_i) \label{RK2}
    \\
    \widehat{y}_{n+1} &= y_n + \Delta t \sum_{i=1}^s \widehat{b}_i f(t_n+c_i\Delta t,u_i), \label{RK3}
\end{align}
\end{subequations} 
where $\Delta t$ is the stepsize, $u_i \in \RR^m$ are the stages which constitute intermediate approximations, and $y_{n+1},\widehat{y}_{n+1}$ are the approximate solutions at $t_{n+1} = t_n + \Delta t$ obtained with $p$th and $\widehat{p}$th order RK methods, respectively,  with $\widehat{p} < p$. 
We refer to the method yielding $y_{n+1}$ as the \emph{advancing method}, and to that yielding
$\widehat{y}_{n+1}$ as the \emph{embedded method}.
Each RK method is characterized by a set of coefficients given in the Butcher tableau
\begin{equation} \label{ButcherTab}
  \begin{tblr}{c|c}
    \cc       & \AAA \\
    \hline
    \,      & \bb^T \\
    \,      & \widehat{\bb}^T
  \end{tblr}, \ \ \ \AAA = [a_{i,j}] \in \RR^{s \times s}, \bb = [b_i], \widehat{\bb} = [\widehat{b}_i], \cc = [c_i] \in \RR^s.
\end{equation}

 It is assumed throughout this work that $c_i = \sum_j a_{i,j}$. If the matrix $\AAA$ in $\eqref{ButcherTab}$ is lower triangular with at least one nonzero diagonal entry, the method is said to be diagonally implicit.
Diagonally implicit Runge-Kutta (DIRK) methods can be unconditionally stable, and their stage equations \eqref{RK1} can be solved sequentially,
whereas in fully implicit Runge-Kutta (IRK) methods, the stage equations are solved simultaneously.  Subclasses of DIRK methods include singly diagonally implicit Runge-Kutta (SDIRK) methods in which the diagonal elements of the matrix $\AAA$ are identical, and explicit first stage diagonally implicit Runge-Kutta (EDIRK) methods; these are DIRK methods satisfying $a_{1,1}=0$. The intersection of these two subclasses is that of explicit first stage singly diagonally implicit Runge-Kutta (ESDIRK) methods. SDIRK methods are an efficient choice when using an algebraic solver in which the linear algebra is handled by a direct method based on matrix factorization, since the factorization can be reused across stages \cite{alexander1977diagonally}.  Meanwhile, EDIRK methods can be designed to have stage order two (cf. section \ref{stgordersec}). We henceforth refer to these subclasses of DIRK methods as DIRK-type methods.

 Another class of IVPs whose solution can be approximated by RK methods is index-1 semi-explicit differential-algebraic equations (DAEs). Such equations take the general form 
\begin{subequations}  \label{DAE}
\begin{align} 
    y' &= f(y,z), \ \ \ y(t_0) = y_0, \\
    0  &= g(y,z),\ \ \ z(t_0) = z_0,
\end{align}
\end{subequations}
for sufficiently smooth functions $f$ and $g$, consistent initial conditions (i.e., $g(y_0, z_0) = 0$), and an invertible Jacobian $\partial g /\partial z$. The unknowns in $\eqref{DAE}$ are the differential variables $y$ and the algebraic variables $z$. The RK iteration for $\eqref{DAE}$ is given by \cite[Sec. VI.1]{wanner1996solving}:
\begin{subequations}\label{RK_DAE}
\begin{align} 
    Y_{n,i} &= y_n + \Delta t \sum_{j=1}^s a_{i,j} f(Y_{n,j}, Z_{n,j}) \\
    0 &= g(Y_{n,j}, Z_{n,i}) \\
    y_{n+1} &= y_{n} + \Delta t \sum_{i=1}^s b_i f(Y_{n,i}, Z_{n,i}) \\
    z_{n+1} &= \left( 1-\sum_{i,j = 1}^s b_i w_{i,j} \right)z_n + \sum_{i,j = 1}^s b_i w_{i,j} Z_{n,j} \label{algebRK1}
\end{align}
\end{subequations}
where $\mathbf{W} := \AAA^{-1} = [w_{i,j}]$. If $\AAA$ is not invertible (e.g., if ESDIRK method is used), then the approximation of the algebraic variable \eqref{algebRK1} can be performed instead by solving the algebraic equation \cite[Sec. VI.1]{wanner1996solving}
\begin{equation}
\label{algebRK2}
    0 = g(y_{n+1}, z_{n+1})
\end{equation}
provided $y_{n+1}$ is available. The DAE $\eqref{DAE}$ can be viewed as a reduced version of the singularly perturbed ODE system 
\begin{subequations}  \label{SPP_general}
\begin{align} 
    y' &= f(y,z) , \ \ \ y(t_0) = y_0 \\
    \varepsilon z'  &= g(y,z), \ \ \ z(t_0) = z_0 
\end{align}
\end{subequations}
but with $\varepsilon=0$. The system \eqref{SPP_general} becomes increasingly stiff as $\varepsilon \to 0$.


Many important IVPs of the form $\eqref{IVP}$ (or $\eqref{SPP_general}$ for small $\varepsilon$) are stiff,
meaning that they cannot be efficiently integrated
with explicit methods due to numerical stability
considerations. This is also the case for DAEs of the form $\eqref{DAE}$.  Efficient numerical approximation of stiff problems and DAEs requires the use of implicit methods, minimally those offering 
unconditional numerical stability. One such stability property is A-stability. As introduced by Dahlquist \cite{dahlquist1963special}, a $k$-step numerical method is A-stable if for every fixed stepsize $\Delta t > 0$, the sequence of approximate solutions $\{y_n\}$ converges to zero as $n \to \infty$ whenever applied to the \emph{linear} test problem 
\begin{equation} \label{testODE}
    y'(t) = \lambda y, \ \ \ \ \ \ \lambda \in \CC \ \ \text{with} \ \ \RE(\lambda) < 0.
\end{equation}

  A fundamental characterization of A-stable methods is that their region of absolute stability include the left complex half-plane. 
  A generalization of the concept of A-stability for $\emph{nonlinear}$ problems is B-stability (see e.g., \cite{wanner1996solving}; Sec. IV.12). However, it has been shown that there exists no B-stable DIRK-type method of order higher than 4 \cite{hairer1980highest}. 
  
 For some classes of problems, A-stability 
 is not sufficient to guarantee a satisfactory numerical solution or may even result in an unstable approximation for certain stiff nonlinear systems \cite{halstead1971mathematical,prothero1974stability}.
 For instance, some initial conditions (e.g., discontinuous initial data for PDEs) may introduce rapid transient phases in the solution. For dissipative problems, fast damping of such transients is desired and may be achieved through L-stability, a concept due to Ehle \cite{ehle1969pade}. When applied to the test problem $\eqref{testODE}$ with a fixed stepsize $\Delta t>0$, an A-stable one-step numerical method is L-stable if at each time instance $t = t_n$, the approximate solution $y_n$ tends to zero as $\RE(\lambda) \to -\infty$, thus preserving the asymptotic behavior of the exact solution to \eqref{testODE}. Another desirable, yet stronger, trait of a numerical method introduced by Prothero and Robinson \cite{prothero1974stability} is stiff accuracy for which the advancing stage $\eqref{RK2}$ is identical to the approximation performed by the last internal stage in $\eqref{RK1}$; in other words, the solution at the next time step $y_{n+1}$ is computed implicitly. Such a property is favorable when approximating the solution to highly stiff, singularly perturbed systems and algebraic variables in DAEs (see for instance \cite{hairer1999stiff} and \cite[Sec. VI., Theorem 1.1]{wanner1996solving}). 

It is known that explicit RK methods cannot be A-stable \cite{nevanlinna1974nonexistence}. Arbitrarily high-order A- and L-stable IRK methods can be constructed based, for instance, on Gauss, Radau and Lobatto quadratures \cite{butcher1964implicit,chipman1971stable,ehle1968high}. Due to their implicit structure, the applicability of these methods may be limited by the computational cost arising from solving the associated coupled nonlinear algebraic systems. Therefore, A- and L-stable DIRK methods may present an efficient alternative. The simplest A-stable DIRK method is the well-known one-stage implicit midpoint rule of order $2$. Various A- and L-stable DIRK-type methods along with their embedded error estimators of up to order 6 have been constructed in the past (e.g., \cite{cooper1979semiexplicit,crouzeix1975approximation,kennedy2016diagonallyRev,kennedy2019diagonally}, and \cite{wanner1996solving}, sec. IV.6). DIRK-type methods of order $6$ include the A-stable stage order $1$ DIRK method (without an error estimator) in 6 stages proposed by Cooper and Sayfy \cite{cooper1979semiexplicit}, and the L-stable stiffly accurate stage order $2$ ESDIRK method in 9 stages with a fifth-order embedding derived in \cite{kennedy2019diagonally}. No sixth-order A-stable SDIRK method exists, nor does there exist an A-stable DIRK-type method of order 7 or higher. Arbitrarily high-order DIRK methods can be constructed using extrapolation, but these methods are only A($\alpha$)-stable (with $\alpha<90^o$) and require a number of stages that grows quadratically with the order \cite[Sec. IV.9]{wanner1996solving}. For instance, an eighth-order DIRK method based on implicit Euler extrapolation must have at least $36$ stages, while the eighth-order A-stable DIRK methods presented in this work have $13-16$ stages.

 This work presents new DIRK-type embedded methods with either fewer stages, new structure, and/or higher order relative to existing methods. The notation we shall adopt in this work to denote each embedded pair is 
 $$\text{TYPE}(s,p)[r]X-[(\widehat{s},\widehat{p})Y]$$
where 
\begin{itemize}
    \item TYPE: the structure of the DIRK-type method, that is DIRK, EDIRK, SDIRK, or ESDIRK.
    \item $s$: number of stages for the advancing method.
    \item $p$: order of convergence for the advancing method.
    \item $r$: stage order of the advancing method.
    \item $X$: stability property of the advancing method; i.e., $A$ for A-stable, $L$ for L-stable, $SA$ for stiffly accurate, or $SAL$ if both stiffly accurate and L-stable.
    \item $\widehat{s}$: number of stages for the embedded error estimator.
    \item $\widehat{p}$: order of convergence for the embedded error estimator.
    \item $Y$: stability property of the embedded error estimator, if any.
\end{itemize}
 Note that in some methods $\widehat{s} = s + 1$ for which we may take $b_{s+1} = 0$ in the embedded pair. Using this notation, the new schemes in this work are:

\begin{outline}
 \1 Sixth-order methods
    \2 \texttt{DIRK(6,6)[1]A-[(7,5)A]}\footnote{While having the same number of stages as the method presented in \cite{cooper1979semiexplicit}, the scheme developed herein is A-stable, double-precision accurate, and has an error estimator, unlike the former where, upon examination, the magnitude of the stability function exceeds 1 for large values of $z$ (i.e., it is not A-stable).} 
    \2 \texttt{DIRK(8,6)[1]SAL-[(8,5)A]} 
    \2\texttt{ESDIRK(8,6)[2]SA-[(8,4)]}\footnote{This method has one fewer stage than the method derived in \cite{kennedy2019diagonally}, but the latter is additionally L-stable.} 
    \2 \texttt{SDIRK(9,6)[1]SAL-[(9,5)A]}

    \1 Seventh-order methods
    \2 \texttt{DIRK(9,7)[1]A-[(9,5)A]} 
    \2 \texttt{DIRK(10,7)[1]SAL-[(10,5)A]} 
    \2\texttt{ESDIRK(10,7)[2]SA-[(10,5)]} 
    \2 \texttt{SDIRK(11,7)[1]SAL-[(11,5)A]}

    \1 Eighth-order methods
    \2 \texttt{DIRK(13,8)[1]A-[(14,6)A]} 
    \2 \texttt{DIRK(15,8)[1]SAL-[(16,6)A]} 
    \2\texttt{ESDIRK(16,8)[2]SAL-[(16,5)]} 
 
\end{outline}

These methods have been found using numerical optimization software, after carefully formulating the order and stability conditions in a way that effectively enables finding their solution by modern optimization algorithms.
We chose to design A-stable and stiffly-accurate methods separately since the need for each of these stability properties is problem-dependent. Moreover, unlike ESDIRK methods, stiffly accurate  DIRK or SDIRK methods may be preferred if direct computation of algebraic variables in DAEs via $\eqref{algebRK1}$ is desired. The coefficients of the new methods are given in the appendix and available electronically\footnote{\url{https://github.com/yousefalamri55/High_Order_DIRK_Methods_Coeffs}}.

These new DIRK-type pairs, which possess combinations of properties not available in any existing
methods, comprise the main contribution of this work.  Additional new contributions herein include:
\begin{itemize}
    \item A new theoretical result regarding the impact of rounding
    errors on A-stability for high-order methods (Section~\ref{sec:stability-rounding});
    \item Reformulations of the problem of finding A-stable RK methods, without which this problem is computationally intractable (Section~\ref{formulation});
    \item New \emph{simplifying conditions} that may facilitate the future construction of very high order RK methods, based on observed properties of the methods found in this work (Section~\ref{SP}).
\end{itemize}

A-stable methods are primarily intended for use in the integration of stiff problems, but DIRK-type methods inevitably suffer from order reduction when applied to stiff problems.  This raises the question of how useful high-order DIRK methods may be.  
A recent study concluded that they can be more efficient than low-order DIRK methods for certain (moderately stiff) problems \cite{kennedy2019diagonally}.  In some cases, a DIRK-type method suffering from order reduction may still be more accurate despite exhibiting a flatter convergence slope (see e.g. \cite[Fig. 5]{2018_wso}).  The methods presented here will enable further investigation into this question, given the especially large difference between their classical order and stage order.

 This paper is organized as follows: in section~\ref{sec2}, we recall the known criteria for A-stability, L-stability, and stiff accuracy. We also recall the conditions for the RK convergence and stage orders. In section~\ref{sec3}, we formulate the optimization problem for finding high-order A-stable and/or stiffly accurate DIRK-type methods.  This problem is, most naturally stated, intractable for modern optimization software, so we describe how the problem can be reformulated in a way that enables finding its accurate numerical solution.  A study of the  stability, structure, and approximation errors for the new schemes as well as strategies for constructing high-order RK methods are presented in section~\ref{sec4}. Several numerical experiments on stiff and differential-algebraic systems are presented in section~\ref{sec5}, and we conclude in section~\ref{sec6} with some remarks and future directions.  

\section{Stability and Order Conditions} \label{sec2}
\subsection{A-stability} 
The application of an $s$-stage RK method of the form $\eqref{RK}$ to the linear test ODE \eqref{testODE} with a stepsize $\Delta t$ reduces to the iteration
\begin{equation}
    y_{n+1} = R(\Delta t \lambda) y_n
\end{equation}
where the rational function $R(\cdot)$ known as the stability function of the RK method. This stability function is given by \cite[Sec. IV.3]{wanner1996solving}:
\begin{align} \label{RKstabfunc}
    R(z) := 1 + z\bb^T(\id - z\AAA)^{-1}\ee = \dfrac{\det(\id - z\AAA + z\ee\bb^T)}{\det(\id-z\AAA)}  =: \dfrac{P(z)}{Q(z)}
\end{align}
for $z \in \CC$, where $\ee = [1,1,\dots,1]^T$ and $\id$ is the $s\times s$ identity matrix. An RK scheme is said to be A-stable if (\cite{wanner1996solving}, Sec. IV.3)
\begin{subequations}  \label{AstabCond}
\begin{align}
    &|R(iy)| \leq 1, \ \ \ \text{for} \ \ \ \ y \in \RR \label{AstabCond_1}\\
    &R(z) \ \ \ \text{is analytic for } \ \ \RE (z) < 0, \ \ z \in \CC. \label{AstabCond_2}
\end{align}
\end{subequations}
Condition $\eqref{AstabCond_1}$ can be reformulated as 
\begin{align} \label{AstabCond2}
     |E(y)| \geq 0, \ \ \text{for\ all} \ y \in \RR
\end{align}
where
\begin{equation*}
    E(y) = |Q(iy)|^2 - |P(iy)|^2
\end{equation*}


A rational function is said to be nondegenerate if its numerator and denominator have no zero in common. Scherer and T\"{u}rke  \cite{scherer1989algebraic} introduced the following equivalent but more practical characterization of A-stability (see also \cite{scherer1994complete}).

\begin{thm}\label{thm:thm1}[Scherer and T\"{u}rke, 1989]
    Let $(\AAA,\bb)$ be the coefficients of an $s$-stage RK method with nondegenerate stability function. The method is A-stable if and only if there exists an $s\times s$ positive-semidefinite symmetric matrix  $\mathbf{R}$ such that $\mathbf{R}\ee = \bb$ and $\MM := \AAA ^T\mathbf{R} + \mathbf{R}\AAA - \bb\bb^T$ is positive-semidefinite. 
\end{thm}

Whereas $\eqref{AstabCond}$ defines A-stability of an RK method in terms of the stability function $\eqref{RKstabfunc}$ and every value $z \in \CC$, the characterization given by Theorem~\eqref{thm:thm1} provides necessary and sufficient conditions for A-stability directly in terms of the coefficients of the RK scheme $(\AAA,\bb)$. Such a characterization is more convenient for designing A-stable schemes via numerical optimization. 

\subsection{A-stability and rounding of coefficients}\label{sec:stability-rounding}
In practical computations, RK methods are used with floating-point systems that are subject to rounding errors.  The effect of rounding errors can be important, for instance, in relation to internal stability for methods with many stages \cite{ketcheson2014internal,kennedy2019diagonally}.  At the same time, the RK coefficients themselves are subject to rounding errors.  This might be especially important with respect to A-stability, since if $|R(z)|\le 1$ is satisfied with equality or near-equality for some values of $z$, then it might be violated in the presence of roundoff.
This issue is especially relevant for methods of very high order,
since the boundary of the stability region tends to get closer and
closer to the imaginary axis as the order is increased \cite{ketcheson2015absolute}.

To investigate effect of rounding RK coefficients, let us introduce perturbed RK coefficients $(\At, \bt)$
with $\|\At-\AAA\|=\OC(\epsilon)$ and $\|\bt-\bb\|=\OC(\epsilon)$,
where $\epsilon$ represents machine precision.  We apply the RK method with coefficients $(\AAA,\bb)$ and the perturbed method with coefficients $(\At,\bt)$ to the test equation \eqref{testODE} and study their difference.
Then the exact and perturbed stages satisfy
\begin{align}
(\id-z\AAA)\UU & = \ee y_n \\
(\id-z\At)\Ut & = \ee y_n, 
\end{align}
respectively, where $\UU=[u_1, \dots, u_s]$.  
A standard result from numerical linear algebra (see e.g.
\cite[Sec.~2.2]{demmel1997applied}) shows that
\begin{align}
    \|\Ut - \UU\| \le \frac{\|(\id-z\AAA)^{-1}\|\cdot\|\At-\AAA\|}{1-\|(\id-z\AAA)^{-1}\|\cdot\|\At-\AAA\|} \|\UU\|.
\end{align}
Next we have
\begin{align}
\yt_{n+1} - y_{n+1} & = z(\bt^T \Ut - \bb^T \UU) \\
        & = z(\bt^T-\bb^T)\UU + z\bt^T(\Ut - \UU) \\
        & = z(\bt^T-\bb^T)(\id-z\AAA)^{-1}\ee y_n + z\bt^T(\Ut - \UU) 
\end{align}
For any DIRK method with all diagonal entries (strictly) positive, it can be shown that both $\|(\id-z\AAA)^{-1}\|$ and $\|z(\id-z\AAA)^{-1}\|$ are
uniformly bounded in the closed left half-plane \cite{calvo2000runge}.  Therefore, for
$\RE (\lambda) \le 0$ in \eqref{testODE}, we have
\begin{align}
    \|\yt_{n+1} - y_{n+1}\| \le \OC(\epsilon) \|y_n\|.
\end{align}
This does not rule out the possibility of a numerical solution
of \eqref{testODE} that grows even when $\RE (\lambda) \le 0$,
but it does guarantee that the growth is extremely slow,
at least for DIRK methods with no vanishing diagonal entries of $\AAA$.

\subsection{L-stability and Stiff Accuracy}  \label{secLstab}
An RK method is L-stable if it is A-stable and its stability function $\eqref{RKstabfunc}$ satisfies (\cite{butcher2016numerical}, Sec. 344):
\begin{equation} \label{LstabCond}
    \lim_{|z| \to \infty} |R(z)| = 0
\end{equation}
A stronger condition that can be used to enforce L-stability is to impose stiff accuracy. An $s$-stage RK method is stiffly accurate if 
\begin{equation} \label{stiffCond}
    a_{s,j} = b_j, \ \ \ \ \ j = 1,2,\dots,s.
\end{equation}
It can be shown that, provided $\AAA$ is nonsingular, $\eqref{stiffCond}$ implies \eqref{LstabCond}; see \cite[Proposition IV.3.8]{wanner1996solving}. Note also that if the RK method is stiffly accurate, then \eqref{RK_DAE} reduces to $y_{n+1} = Y_{n,i}$ and $z_{n+1} = Z_{n,i}$.

\subsection{Order Conditions}
The order conditions for an RK method are a set of algebraic relations among its coefficients constituting the necessary and sufficient conditions for an RK method to yield a given order of convergence. As the order of an RK method increases, the number of order conditions increases rapidly. Despite the known algebraic structure of the RK order conditions \cite{butcher1972algebraic,hairer1974butcher}, and the existence of many simplifying assumptions for these order conditions, the analytic design of high-order RK schemes is still extremely challenging. Imposing A- or L-stability renders it even
more difficult.

An RK method with coefficients $(\AAA,\bb)$ is of order $p$ if and only if \cite[Sec. II.2]{hairer1993solving}
\begin{align}\label{OCs}
{\boldsymbol\tau}^{(q)}(\mathbf{A},\mathbf{b}) &:=  \left( \sum_{j=1}^s b_j\Phi_{i,j}^{(q)} - \dfrac{1}{\gamma(t)} \right) = 0 , \ \ \text{for} \ \ t \in T_q, \ q = 1,2,\dots,p.
\end{align}
where $T_q$ is the set of all rooted trees of order $q$, $b_j$ are the RK advancing weights in $\eqref{RK2}$, $\gamma(t)$ is the density of the tree $t$, and $\Phi_{i,j}^{(q)}$ are certain nonlinear products of the RK coefficients $a_{i,j}$. For a given order $q$, the number of the required conditions is precisely the cardinality of $T_q$, denoted by $|T_q|$. Hence, the total number of order conditions for a $p$th order RK scheme is the cumulative sum of the cardinality of each tree of order $1,\dots,p$. In Table~\ref{tab:OCTable} we recall the cardinality of the rooted trees and the total number of required order conditions for an RK method to be of order $p$.  Notice that the maximal degree of the polynomial expressions appearing in the order conditions is also equal to the order.

 \subsection{Stage Order} \label{stgordersec}
 Besides the classical order (that is, the order of convergence) of an RK method $p$, the $\emph{stage order}$, call it $r$, plays an important role in the behavior of the global error when approximating the solution of stiff ODEs or DAEs. Let $(\AAA,\bb)$ be the coefficients of a $p$th order RK method and suppose $\ell$ is the largest positive integer such that 
\begin{equation} \label{Cm}
    \mathcal{C}_k(\AAA) := 
    \AAA \cc^{k-1}- \dfrac{\cc^k}{k} =0, \ \ \ \ \ \ k = 1,2,\dots,\ell.
\end{equation}
holds, where the power $k$ here is assumed to be componentwise. Then, the stage order of the RK method is \cite[Sec. IV.15]{wanner1996solving}:
\begin{equation} \label{stgOrd}
    r = \min\left\{p,\ell\right\}
\end{equation}
All RK methods used in practice have stage order at least 1
(i.e., they satisfy $c_i = \sum_j a_{i,j}$). It is known that the stage order of an SDIRK (or DIRK with nonzero diagonal elements) method is at most $1$, and that of an ESDIRK or EDIRK method is at most $2$ \cite[Sec. IV.15, Ex. 1]{wanner1996solving}.

\begin{table}
\centering
    \caption{The cardinality of a rooted tree $T_q$, denoted by $|T_q|$, and the cumulative number of order conditions for a $p$th order RK scheme.} 
    \label{tab:OCTable} 
\begin{tabular}{c|cccccccccc} \toprule 
    {$q$} & {$1$} & {$2$} & {$3$} & {$4$} & {$5$} & {$6$} & {$7$} & {$8$} & {$9$} & {$10$}\\ \midrule
    {$|T_q|$} & {$1$} & {$1$} & {$2$} & {$4$}  & {$9$} & {$20$} & {$48$} & {$115$} & {$286$} & {$719$}\\ \midrule
    {$\sum_{q=1}^p |T_q|$} & {$1$} & {$2$}  & {$4$} & {$8$} & $17$ & $37$ & $85$ & $200$ & $486$ & $1205$
       \\ \bottomrule
\end{tabular}
\end{table}



\section{Optimization Approach} \label{sec3}
In this section, we describe the optimization approach used to find the coefficients of the new DIRK pairs. Henceforth, let $(\AAA,\bb[\widehat{\bb}])$ be the desired RK coefficients, ${\boldsymbol\tau}^{(q)}(\AAA,\bb[\widehat{\bb}])$ be the set of RK $q$th order conditions given by $\eqref{OCs}$, and ${\boldsymbol{\mathcal{C}}}_k(\AAA)$ be the $k$th stage order condition given by \eqref{Cm}.


\subsection{Problem formulation} \label{formulation}
The optimization problem for finding a $p$th order, stage order $r$, A-stable, and stiffly accurate (thus L-stable if $\AAA$ is chosen such that it is nonsingular) DIRK-type scheme with $s$ stages is
\begin{equation} \label{OpProb-naive}
\begin{aligned} 
\text{Choose } & \RRR,\AAA \in \RR^{s \times s}, \bb [\widehat{\bb}] \in \RR^s\\
\text{s.t.} \ \ \  
&  {\boldsymbol\tau}^{(q)}_l(\AAA,\bb) = 0 &   (\text{for} \ \ t \in T_q, \ l = 1,\dots, |T_q|, \ \ q = 1,\dots,p)\\
& {\boldsymbol{\mathcal{C}}}_k(\AAA) = 0 &  (\text{for} \ \ k = 1,\dots, r)\\
& |1 + z\bb^T(\id - z\AAA)^{-1}\ee|\le 1 & \forall z\in \CC^- \\
  &a_{i,j} = 0 &  (\text{for} 
  \ j > i) \\
  & \sigma a_{s,j} =  \sigma b_j &  (\text{for} 
  \ j = 1,2,\dots,s).
\end{aligned}
\end{equation}
Here the constraints are, in order, the classical order conditions, stage order conditions, A-stability (with $\CC^-$ denoting the closed left half of the complex plane), DIRK structure, and stiff accuracy (cf. section~\ref{secLstab}).  We set $\sigma=1$ to impose stiff
accuracy or $\sigma=0$ to omit this condition.
The (E)(S)DIRK structure can be imposed additionally if desired. As mentioned in section~\ref{stgordersec}, for EDIRK or ESDIRK methods, $r \in \{1,2\}$, while for SDIRK (or DIRK with nonzero diagonal elements) methods $r$ cannot exceed $1$. 

The problem in this most natural form cannot be solved by numerical
optimization, due to the A-stability constraint which must be imposed
over an uncountable set.  One can reduce this set from the whole left
half-plane to just the imaginary axis, but this does not essentially
improve matters.  Instead, we can apply Theorem~\ref{thm:thm1} in order
to replace this constraint by a (finite) semi-definite constraint. This yields:
\begin{equation} \label{OpProb}
\begin{aligned} 
\text{Choose } & \RRR,\AAA \in \RR^{s \times s}, \bb [\widehat{\bb}] \in \RR^s\\
\text{s.t.} \ \ \  
&  {\boldsymbol\tau}^{(q)}_l(\AAA,\bb) = 0 &   (\text{for} \ \ t \in T_q, \ l = 1,\dots, |T_q|, \ \ q = 1,\dots,p)\\
& {\boldsymbol{\mathcal{C}}}_k(\AAA) = 0 &  (\text{for} \ \ k = 1,\dots, r)\\
& \AAA ^T\mathbf{R} + \mathbf{R}\AAA - \bb\bb^T \succcurlyeq 0 \\
& \mathbf{R}\ee = \bb \\
& \RRR = \RRR^T \\
& \RRR \succcurlyeq 0 \\
  &a_{i,j} = 0 &  (\text{for} 
  \ j > i) \\
  & \sigma a_{s,j} =  \sigma b_j &  (\text{for} 
  \ j = 1,2,\dots,s).
\end{aligned}
\end{equation}
Now the problem has a finite set of constraints, but it is not yet
tractable for numerical solvers.
Since problem \eqref{OpProb} involves nonlinear and semidefinite constraints, it is natural to apply either semidefinite programming or nonlinear programming algorithms for finding its solution.
But existing semidefinite programming tools are not designed to handle non-convex constraints, while many nonlinear programming tools do not accept semidefinite constraints, so most optimization codes will not accept
this problem.  In our experience, those solvers that can accept
the problem in this form (which treat the constraints as a black box)
consistently fail to find any solution for $p\ge 6$.  This is
perhaps unsurprising due to the complexity and quantity of the constraints.

To make the problem suitable for solvers that employ more information
(e.g. gradients of the constraints), we further
reformulate the semidefinite matrix constraints
via Cholesky decomposition.  Recall that a matrix $\pmb{X} \in \RR^{n\times n}$ is symmetric positive semidefinite if and only if there exists a lower triangular matrix $\pmb{L} \in \RR^{n\times n}$ whose diagonal elements are nonnegative and $\pmb{X} = \pmb{L}\pmb{L}^T$.  This converts the semidefinite constraints to polynomial constraints, at the cost of introducing more decision variables.  Yet even in this form, we have found that numerical optimization tools are unsuccessful in finding solutions.
We have therefore reformulated the problem once more by relaxing
constraints through the introduction of slack variables,
resulting in the following (final) formulation:
\begin{equation} \label{OpProbReform}
\begin{aligned}
\min_{\RRR,\AAA,\bb,\LM,\LR} \quad & \sum_{q=1}^p \sum_{l=1}^{|T_q|}  |\lambda^{(q)}_{l}| + \sum_{k=1}^r \sum_{j=1}^{s}  |\lambda^{(k)}_{j}| +\sum_{i=1}^s \sum_{j=1}^s |\lambda_{i,j}^M| + \sum_{i=1}^s |\lambda_{i}^{Re}|  \\ 
\quad & + \sum_{i=1}^s \sum_{j=1}^s |\lambda_{i,j}^R| + \sum_{i=1}^s  \lambda_{i}^{R} + \sum_{i=1}^s \lambda_{i}^{M} + \sum_{i=1}^{s(s-1)/2} |\lambda_{i}^{a}| \\
\textrm{s.t.} \quad & {\boldsymbol\tau}^{(q)}_l(\AAA,\bb) = \lambda^{(q)}_l \ \ \ \ \ \  \ \ \ \ \ \  (\text{for} \ \ t \in T_q, \ l = 1,\dots, |T_q|, \ \ q = 1,\dots,p)\\
    \quad & {\boldsymbol{\mathcal{C}}}_k(\AAA) = \lambda^{(k)}_j \ \ \ \ \ \ \ \ \ \ \ \ \ \ \ \ \ \ \ \ \ \ \ \ \ \  \ \ \ \ \  (\text{for} \ \ k = 1,\dots,r, \ j = 1,\dots, s)\\
  \quad & \MM - \LM \LM^T = \pmb\Lambda^M \ \ \ \ \ \ \ \ \ \ \ \ \ \ \ \ \ \ \ \ \ \ \ \ \ \ \ \ \ \ \ \ \ \ \ \ \ (\pmb\Lambda^M = [\lambda_{i,j}^M] \in \RR^{s \times s})\\
  \quad & \RRR \ee - \bb = \pmb\Lambda^{Re}\ \ \ \ \ \ \ \ \ \ \ \ \ \ \ \ \ \ \ \ \ \ \ \ \ \ \ \ \ \ \ \ \ \ \ \ \ \ \ \ \ \ \ \ \ ( \pmb\Lambda^{Re} = [\lambda_i^{Re}] \in \RR^s)\\
  \quad & \RRR - \LR \LR^T = \pmb\Lambda^R \ \ \ \ \ \ \ \ \ \ \ \ \ \ \ \ \ \ \ \ \ \ \ \ \ \ \ \ \ \ \ \ \ \ \ \ \ \ \ \ 
  \ (\pmb\Lambda^R = [\lambda_{i,j}^R] \in \RR^{s \times s}) \\ 
  \quad & \RRR = \RRR^T \\
  \quad & \diag(\LM) \geq \pmb\Lambda^{M} \ \ \ \ \ \ \ \ \ \ \ \ \ \ \ \ \ \ \ \ \ \ \ \ \ \ \ \ \ \ ( \pmb\Lambda^{M} = [-\lambda_i^{M}] \in \RR^s, \lambda_i^{M} \geq 0) \\
  \quad & \diag(\LR) \geq \pmb\Lambda^{R} \ \ \ \ \ \ \ \ \ \ \ \ \ \ \ \ \ \ \ \ \ \ \ \ \ \ \ \ \ \ \ \ \ \  ( \pmb\Lambda^{R} = [-\lambda_i^{R}] \in \RR^s, \lambda_i^{R} \geq 0) \\
  \quad & a_{i,j} = \lambda_{i}^{a} \ \ \ \ \ \ \ \ \ \ \ \ \ \ \ \ \ \ \ \ \ \ \ \ \ \ \ \ \ \ \ \ \ \ \ \ \ \ \ \ \ \ \ \ \  \ \ \ \ \ \ \ \ \ \ \ \ \ \ \ \  \ \ \ \ (\text{for}
  \ j > i)\\
\quad & \sigma a_{s,j} = \sigma b_j  \ \ \ \ \ \ \ \ \ \  \ \ \ \ \ \ \ \ \ \ \ \ \ \ \ \ \ \ \ \ \ \ \ \ \ \ \ \ \ \ \ \ \ \ \  \ \ \ \ \ (\text{for} 
  \ j = 1,2,\dots,s).
\end{aligned}
\end{equation}

where $\MM := \AAA ^T\mathbf{R} + \mathbf{R}\AAA - \bb\bb^T$, $p$ is the order of the desired scheme, $r$ is its stage order, and $s$ is the number of stages. A solution of problem $\eqref{OpProbReform}$ represents a solution of problem \eqref{OpProb} if the sum of the introduced slack variables (i.e., the objective function) is equal to zero.

To find an embedded error estimator for a given scheme, we simply fix
the matrix $\AAA$ and then solve \eqref{OpProbReform} (with $\sigma=0$ and the desired value of the order $\hat{p}$ for the error estimator).
This problem is considerably simpler since the order conditions are all linear when $\AAA$ is fixed.
For very high order methods, this latter problem may be so constrained that there are no solutions
other than the original method, in which case we may reduce the desired order or increase the number of
stages in order to find a non-defective error estimator.
Since the stability function of an ESDIRK methods is degenerate, designing A-stable embedding for them using Theorem~\ref{thm:thm1} is not possible. However, when obtained using the formulation delineated above, the embedded error estimators for the ESDIRK methods possess practically large stability regions.  
 

Our goal is to find methods that satisfy the desired constraints and have the fewest stages possible.
In Table \ref{tbl:stages} we list the best known bounds on the minimum number of stages required to construct a DIRK method of a given order.  Up to order six, the bounds are known to be tight as they are achieved by existing methods in the literature.  For order seven and higher, the bounds are based on a result related to rational functions with real poles, which states that the degree of approximation to the exponential is at most one greater than the degree of the numerator \cite{norsett1977attainable}. These bounds are valid when considering DIRK schemes without requiring A-stability.  For construction of A-stable DIRK-type schemes, there are no available bounds; indeed, there is no available proof that A-stable DIRK-type methods of order seven (or higher) exist.
\begin{table}[!h!]
\centering
    \caption{Known lower bounds on the minimum number of stages required to construct a DIRK method of order $p$.  Numbers in parentheses are tighter bounds known to be required for A-stable schemes.}
    \label{tbl:stages} 
\begin{tabular}{c|cccccccccc} \toprule 
    $p$            & 1 & 2 & 3 & 4 & 5    & 6 & 7 & 8 & 9 & 10\\ \midrule
    $s_{\rm{min}}$ & 1 & 1 & 2 & 3 & 4(5) & 5(6) & 6 & 7 & 8 & 9 \\ \midrule
\end{tabular}
\end{table}

\subsection{Optimization Strategy}
A combination of global and local nonlinear solvers was used to solve $\eqref{OpProbReform}$ numerically. Global optimization has been employed previously in the context of designing RK schemes (see for instance \cite{ketcheson2009optimal,ruuth2006global,ruuth2004high}). Due to the increase in the number and intricacy of the imposed constraints herein, as well as the commensurate computational complexity introduced by the nature of global optimization algorithms, a guarantee of global optimality becomes rather difficult for high-order RK schemes with additional stability requirements. We found it useful to employ a two-phase search strategy, first using a global solver to find an approximate solution and then using a local optimizer to improve the solution accuracy.

We attempted to find methods satisfying the desired constraints using the smallest possible number of stages.  The lower bounds for the number of stages for DIRK methods given in table~\ref{tbl:stages} were used as a starting point to search for solutions.  If no solution could be found,
the number of stages was incremented and a new search
conducted until an RK method satisfying the constraints 
(to within double precision) was obtained.

To solve \eqref{OpProbReform}, \texttt{BARON} \cite{sahinidis1996baron}, a
deterministic, global optimization solver that employs branch-and-bound type algorithms, through its \texttt{GAMS} \cite{bussieck2004general} interface was used in the first phase to obtain an approximate solution. In essence, this algorithm solves a relaxed version of the primal problem over successively refined partitions of the search space and determines an upper and lower bound on the primal objective. As the refinement of the search domain continues, the gap between the upper and lower bounds on the objective decreases until it meets the termination criterion. In \texttt{BARON}, both the absolute \texttt{EpsA} and relative \texttt{EpsR} termination tolerances, which control the desired gap between the upper and lower bounds on the global solution,  were set to be $10^{-12}$. The feasibility tolerance \texttt{AbsConFeasTol}, the allowable tolerance for the constraint violation, was set to $10^{-10}$. In addition to these specified parameters, an ample maximum CPU time limit was enforced as termination criterion. As subsolvers in \texttt{BARON}, the nonlinear programming solver was chosen to be \texttt{CONOPT} and the selected linear programming solver was \texttt{CPLEX}.

In addition to the aforementioned specifications, we provide some remarks that may be helpful regarding the use of global solvers based on our experience for solving \eqref{OpProbReform} in $\texttt{BARON}$.  In general, there is a tradeoff between using equality constraints (or decision variables) and eliminating as many as possible of them; this can make the problem either more or less tractable to a given solver, and our experience in this respect shows that this is highly problem-, constraint-, and solver-dependent and is extremely difficult to predict. First, although it reduces the number of the order conditions, imposing the simplifying assumptions (see e.g.,
\cite[Sec. 321]{butcher2016numerical}) for the RK order conditions (besides the stage order condition \eqref{Cm}) was found to hinder the global search. Second, it was observed that supplying feasible bounds for the decision variables increases the efficiency of the solver as it, expectedly, narrows the search domain. Third, when searching for sixth and seventh order DIRK-type schemes, the solution to $\eqref{OpProbReform}$ was obtained while completely eliminating the upper triangular elements of the matrix $\AAA$ from the optimization problem (rather than adding additional constraints that equate them to 0 to enforce the DIRK structure). A similar elimination was done for the stiff accuracy constraint. On the other hand, finding the eighth order methods with the global solver was only possible while retaining these variables and using explicit equality constraints along with the associated slack variables. Otherwise, the global search would not terminate. This may suggest the need for the additional relaxation introduced by the slack variables when using global solvers to solve problems with a large number of constraints and few degrees of freedom.  Finally, while $\texttt{BARON}$ has the capability for choosing initial guesses judiciously without a user input, the use of random initial guesses (consisting of a random lower-triangular matrix $\AAA$ and random vector $\bb$) made it possible to find different sets of solutions, which was necessary since some of the solutions obtained in the first phase (global search) could not be made sufficiently accurate in the second phase (local search) or they had abnormally large or small coefficients (in which case they would not be good candidates for RK methods as solving for the stages would become extremely difficult or computing $\AAA^{-1}$ would be prone to inaccuracy).  Furthermore, due to the difficulty of finding individual solutions, it was not possible to add additional constraints that, for instance, minimize some measures of the leading truncation error for the RK methods.

Once an approximate solution was obtained via \texttt{BARON}, a sequence of different local nonlinear solvers (including \texttt{CONOPT} \cite{drud1994conopt}, \texttt{IPOPT} \cite{biegler2009large}, \texttt{KNITRO} \cite{byrd2006k}, and \texttt{SNOPT} \cite{gill2005snopt}) were used to further reduce the constraint residuals as far as possible.  In $\texttt{IPOPT}$, relaxation of fixed variables was activated by setting the $\texttt{fixed\_variable\_treatment}$ parameter to $\texttt{relax\_bounds}$.  For $\texttt{KNITRO}$, the chosen algorithm was the sequential quadratic programming (SQP) algorithm. In all local solvers, the desired relevant tolerances were set to $10^{-20}$ (although not always achieved by the solvers, they would return the best possible solutions in terms of accuracy). Finally, the obtained coefficients were further refined via local search in $\texttt{Mathematica}$ using the function $\texttt{FindMinimum}$ with respect to the RK order conditions until the maximum residual of the order condition is of order $10^{-16}$ or less.  In this last step, although the A-stability constraint was not explicitly imposed, only the last four to five digits (out of sixteen digits) of the coefficients change and the stability of each resulting method had to be carefully verified afterward.

\section{New High-order DIRK Methods} \label{sec4}
We present eleven new DIRK-type schemes of orders 6-8 with their embedded error estimators.  All of them are A-stable; some of them are additionally L-stable and/or  stiffly accurate. The schemes of 7th and 8th orders and all of the SDIRK schemes are the first of their kind. The ESDIRK schemes have stage order 2. The abessecia are computed  as $c_i = \sum_j a_{i,j}$ and, for the DIRK methods, satisfy $0 \leq c_i \leq 1$, for $i = 1,\dots,s$.  The coefficients of the embedded pairs up to 16 digits are included in the appendix.



\subsection{Stability}

\subsubsection{Absolute Stability}
\begin{figure}[!ht!]
\includegraphics[width=6.5in]{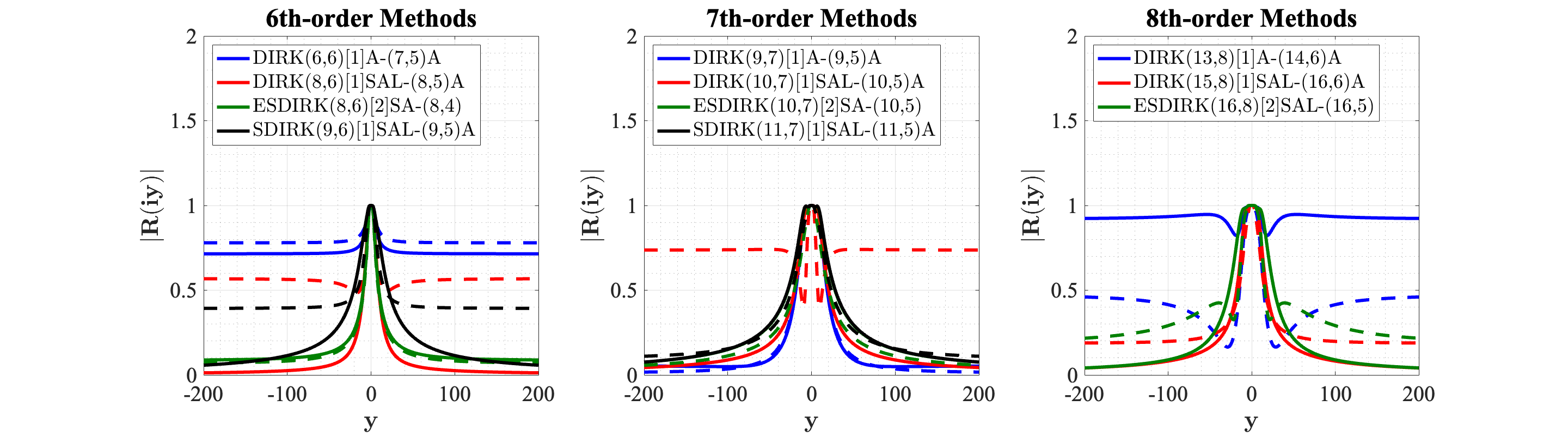}
\caption{The modulus of the stability function \eqref{RKstabfunc} for 6th- (left), 7th- (middle), and 8th- (right) order schemes. The solid curves are for the advancing method and the dashed curves is for the corresponding embedded error estimator.}
\label{fig:StabFunc}
\end{figure}

The schemes found are guaranteed to satisfy the
semidefinite constraints that are equivalent to A-stability, but only up to a certain numerical precision (typically 15-16 digits).  For all schemes, we have further verified their
A-stability, by calculating the magnitude of the stability function on the imaginary axis.
As shown in Figure~\ref{fig:StabFunc}, the modulus of the stability function of each of the A-stable methods (solid curves for the advancing method, dashed curves for the embedded error estimators) is less than or equal to unity in magnitude. Additionally, the first 
 column of Table~\ref{tbl:StabFuncVals} lists the magnitude of the stability function at $z=\pm \infty$ for each method and its embedded error estimator, confirming L-stability of all methods for which condition $\eqref{LstabCond}$ was imposed. None of the error estimators for the ESDIRK methods is A-stable (cf. section \ref{formulation} for explanation), but these estimators do have large regions of absolute stability as $|R(iy)| 
 > 1$ for $|y|> 10^9$,  $|y|> 10^{10}$, and $|y|> 10^6$ for the embedded error estimators of the 6th, 7th, and 8th-order ESDIRK schemes, respectively. The embedded error estimators for the rest of the methods are A-stable.

   \begin{table}[!h!]
   \small
    \centering
\caption{The stability function \eqref{RKstabfunc} values at $\pm \infty$ and maximum values for the internal stability functions \eqref{internal_1} and \eqref{internal_2}.}
\label{tbl:StabFuncVals}
    \begin{tblr}{ vline{2}={\heavyrulewidth},
                 colspec={Q[l, mode=math] ccc}}
    \toprule
    \text{\ \ \ \ \ \ \ \ \ \ \ \ \  Method}           & $|R(\pm \infty)|$ & $\max_{j,y} |\rho_j(iy)|$ & $\max_{j,y} |\theta_j(iy)|$\\ 
    \midrule
   \texttt{DIRK(6,6)[1]A-[(7,5)A]}    &  0.71 [0.78]    & 1.10 [1.10] &  0.40 [0.40] \\
    \texttt{DIRK(8,6)[1]SAL-[(8,5)A]}     &   0.00 [0.57]   & 1.08    &  0.31 [0.31]  \\
    \texttt{ESDIRK(8,6)[2]SA-[(8,4)]}     &  4.77   & 2.33   & 0.42 [0.41]  \\
   \texttt{SDIRK(9,6)[1]SAL-[(9,5)A]}     &  0.00 [0.39]    & 1.29   &  0.81 [1.00]    \\
    \texttt{DIRK(9,7)[1]A-[(9,5)A]}     &    0.06 [0.01]     & 1.11  &  1.19 [1.16]   \\
    \texttt{DIRK(10,7)[1]SAL-[(10,5)A]}      &  0.00 [0.74]    & 1.23  &  0.92 [0.39]   \\
    \texttt{ESDIRK(10,7)[2]SA-[(10,5)]}    &  0.01       & 11.27    &  0.37 [0.39]  \\
    \texttt{SDIRK(11,7)[1]SAL-[(11,5)A]}    &  0.00 [0.09]   & 1.02   &   0.70 [0.63]    \\
    \texttt{DIRK(13,8)[1]A-[(14,6)A]}      &  0.92 [0.48]     & 2.60 [2.59] &  0.71 [0.62]   \\
    \texttt{DIRK(15,8)[1]SAL-[(16,6)A]}     &   0.00 [0.19]   & 4.95 [4.95] &  0.51 [0.35]    \\
    \texttt{ESDIRK(16,8)[2]SAL-[(16,5)]}     &   0.00    & 12.52    &  0.34 [0.33]   \\
    \bottomrule
    \end{tblr}
    \end{table}

\subsubsection{Internal Stability}  
Given the large number of stages used by the presented schemes, internal stability is an important concern.  Internal stability governs the amplification of perturbations within the stages of a single step of the RK method.  The internal stability functions \cite{kennedy2019diagonally} 
\begin{align} \label{internal_1}
    \rho_j(z) & = [(\id-z\AAA)^{-1}\ee]_j
\end{align}
indicate how a perturbation to $y_n$ is amplified up to stage $u_j$, while the functions \cite{ketcheson2014internal}
\begin{align}  \label{internal_2}
    \theta_j(z) & = [\bb^T(\id-z\AAA)^{-1}]_j
\end{align}
indicate how a perturbation to $u_j$ affects the next step
value $y_{n+1}$.  It is desirable to ensure that these functions do not take large values in the left half-plane, and particularly as $|z|\to\infty$.
For DIRK and SDIRK methods, both $\rho_j$ and $\theta_j$ are guaranteed to vanish at infinity as $a_{ii}\ne 0$ for each $i$.  The maximum values of these internal stability functions are listed for each of the methods in Table \ref{tbl:StabFuncVals}.

\subsection{Error Coefficients}
One approach to quantify the error for a $p$th order RK method is to examine the norm of the leading order truncation error coefficients. These coefficients are the norm of the $(p+1)$th and $(p+2)$th order terms of the left hand side of the order conditions $\eqref{OCs}$:

\begin{align} \label{E1}
    E_{\infty}^{(p+1)} := \| {\boldsymbol\tau}^{(p + 1)}(\AAA,\bb[\widehat{\bb}] )  \|_{\infty}
\end{align}
and 
\begin{align}\label{E2}
    E_{\infty}^{(p+2)} := \| {\boldsymbol\tau}^{(p + 2)}(\AAA,\bb[\widehat{\bb}])  \|_{\infty}
\end{align}
 Another important measure in the context of RK methods is the largest (in magnitude) coefficient in the Butcher tableau \eqref{ButcherTab} given by
\begin{equation}\label{D}
    D := \max\{|a_{i,j}|,|b_i [\widehat{b}_i]|,|c_i|\}
\end{equation}
 These quantities are provided in Table~{\ref{tbl:Errs}}.


 \begin{table}
 \small
    \centering
\caption{The error measures \eqref{E1}-\eqref{E2} in the $L^{\infty}$ norms and the maximum coefficients \eqref{D} associated the DIRK-type schemes and their embedded error estimators.}
\label{tbl:Errs}
    \begin{tblr}{vline{2}={\heavyrulewidth},
                 colspec={Q[l, mode=math] lllll},
                 row{1} = {mode=math}}
    \toprule
    \text{\ \ \  \ \ \ \  \  \ \ \  Method}             & \ \ \  \  \  \ \ \ E_{\infty}^{(p+1)}  & \ \ \  \  \  \ \ \ E_{\infty}^{(p+2)} & \ \ \  \  \  \ \ \ D \\ 
    \midrule
    \texttt{DIRK(6,6)[1]A-[(7,5)A]}    &  1.75e-3 [9.19e-4]    & 5.16e-3 [1.96e-3] &  1.00e+0 [1.00e+0] \\
    \texttt{DIRK(8,6)[1]SAL-[(8,5)A]}     &    3.83e-4 [7.03e-4]   & 9.99e-4 [1.09e-3]  & 1.00e+0 [1.00e+0]   \\
    \texttt{ESDIRK(8,6)[2]SA-[(8,4)]}     &  1.07e-3 [3.94e-4]   & 1.92e-3 [8.00e-4]  &  1.21e+0 [1.21e+0]  \\
   \texttt{SDIRK(9,6)[1]SAL-[(9,5)A]}     &  1.84e-4 [9.28e-4]    & 2.42e-4 [8.03e-4] &  1.00e+0 [1.00e+0]   \\
    \texttt{DIRK(9,7)[1]A-[(9,5)A]}     &    6.55e-5 [3.26e-5]     & 4.83e-5 [1.90e-5] &  1.19e+0 [1.16e+0]   \\
    \texttt{DIRK(10,7)[1]SAL-[(10,5)A]}      & 1.96e-5 [3.68e-4]    & 4.17e-5 [5.92e-4] &  1.00e+0 [1.00e+0]   \\
    \texttt{ESDIRK(10,7)[2]SA-[(10,5)]}    &   6.64e-5 [3.26e-4]     & 1.04e-4 [4.91e-4] & 1.00e+0 [1.00e+0]   \\
    \texttt{SDIRK(11,7)[1]SAL-[(11,5)A]}    &  1.29e-5 [7.13e-5]   & 2.86e-5 [9.43e-5] & 1.03e+0 [1.03e+0]      \\
    \texttt{DIRK(13,8)[1]A-[(14,6)A]}      &  8.99e-5 [1.30e-4]     & 9.60e-5 [2.44e-4] &  1.00e+0 [1.00e+0]   \\
    \texttt{DIRK(15,8)[1]SAL-[(16,6)A]}     &   6.08e-5 [1.81e-4]  & 1.01e-4 [3.87e-4] &   1.00e+0 [1.00e+0]   \\
    \texttt{ESDIRK(16,8)[2]SAL-[(16,5)]}     &  3.12e-6 [6.82e-5]    & 3.67e-6 [7.00e-5] & 1.00e+0 [1.00e+0]    \\
    \bottomrule
    \end{tblr}
    \end{table}

\subsection{Structural Properties} \label{SP}
In the construction of high-order RK methods, it is common to make use of
\emph{simplifying assumptions}.  These are conditions on the coefficients that may not
be necessary but make it possible to satisfy many order conditions simultaneously
with few degrees of freedom.  In the present work, we have not made use of any simplifying
assumptions, since our goal is to find methods with the fewest stages possible by imposing
only the necessary conditions for a given order of accuracy.  In this section we study
relationships that are in fact satisfied (to within rounding errors) by the coefficients of the methods we have found,
but which are not (at least explicitly) necessary.  These relationships provide some insight
into structural conditions (similar to but different from the simplifying assumptions in the
literature) that may facilitate the construction of high-order RK schemes in the future, whether
by manual solution of the order conditions or by numerical optimization.

The discussion in this section is greatly simplified by considering a form of the order
conditions different from that in which they are usually written.  By taking certain linear
combinations of the conditions for order $p$, one can construct an equivalent set of conditions 
that has the form
\begin{align}
    \sum_i b_i c_i^k & = \frac{1}{k+1} & k = 0, 1, \dots, p \\
    \sum_i b_i v_i(\tau) & = 0 & \mbox{ for each } \tau \in T_p.
\end{align}
Here each vector $v(\tau)$ is a function of the RK matrix $\AAA$ only, and $T_p$
is the set of all rooted trees of order at most $p$.
A straightforward construction of the order conditions in this form is given in
\cite[Section~4.2]{albrecht1996}.  These order conditions are most easily expressed
using the \emph{stage order residual vectors}
\begin{align}
    \gamma^{(k)} & = k \AAA \cc^{k-1} - \cc^k & k \ge 1.
\end{align}
One sees that if $\gamma^{(k)}_i=0$ for $1\le k \le q$ then the stage $u_i$
is an approximation to $y(t_n+c_i\Delta t)$ with local order of accuracy $q$,
and we say that stage $i$ has order $q$.
To simplify the discussion in the rest of this section, when we say that the method
coefficients satisfy an exact property, we mean that it is satisfied up to
rounding errors.

\subsubsection{Stage order 2 and zero weights}
It turns out that for most of the methods found herein (all except the \texttt{DIRK(6,6)[1]A} method), 
the following property holds:
\begin{align} \label{only_use_so2}
    \gamma^{(2)}_i \ne 0 \iff b_i = 0.
\end{align}
This means that any condition of the form
\begin{align} \label{so2_simp}
    \sum_i b_i w_i \gamma^{(2)}_i & = 0
\end{align}
is automatically satisfied for any vector $w$.  For instance, up to order
8 there are 37 (out of 200) conditions that have this form.  Most of the
methods found in this work have $b_1=0$ and $\gamma^{(2)}_i=0$ for $i>1$.
This is a very common structural condition that is also satisfied for
instance by many explicit 5th-order methods.  For method construction,
if one imposes \eqref{only_use_so2} \emph{a priori}, then the corresponding 
order conditions \eqref{so2_simp} can be neglected.

\subsubsection{Structure of \texttt{DIRK(15,8)[1]SA}}
Note that condition \eqref{only_use_so2} means that only stages with order
of accuracy two or higher are used in the construction of the new solution.
Even more simplification can be gained if all stages used in the solution
update have order higher than two.  This turns out to be the case for the
\texttt{DIRK(15,8)[1]SA} method, which has additional structure that we discuss here.

Let $I=\{1,2,9\}$. For this method, $b_j=0$ for $j\in I$, while
$\gamma^{(2)}_j=\gamma^{(3)}_j = 0$ for all $j \notin I$.  This means
that \eqref{so2_simp} holds for all vectors $w$, and additionally that
(for any $w$) we have
\begin{align} \label{so3_simp}
    \sum_i b_i w_i \gamma^{(3)}_i & = 0.
\end{align}
Among the 200 conditions for order eight, there are 17 conditions that take this form.

Additionally, for the \texttt{DIRK(15,8)[1]SA} method we find that $(\AAA\gamma^{(2)})_j=0$
for all $j\notin I$.  This means further that for any $w$ we have
\begin{align} \label{Ag2_simp}
    \sum_{i,j} b_i (a_{ij} \gamma^{(2)}_j) w_i & = 0.
\end{align}
Another 17 of the 200 conditions for order eight have this form.

\subsubsection{Confluence of low-order stages}
An additional kind of simplifying structure was observed in a method developed in this work but which was not finally selected for presentation.  In that method, which is an 8th-order A-stable DIRK method with 13 stages, we observed that $\gamma^{(2)}_j=0$ for $j\notin \{1,9\}=:I$.  As in the case of \texttt{DIRK(15,8)[1]SA}, for this method $b_j=0$ for
$j\in I$.  Additionally, for this method $c_1=c_9$.  Therefore, $\gamma^{(2)}$ is an eigenvector of the matrix
$C:=\diag(c)$.  Thus order conditions involving factors of
the form $C^k \gamma^{(2)}$ are automatically satisfied if the corresponding condition involving $C \gamma^{(2)}$ is
satisfied.  This structure is also seen in the high-order explicit methods of Zhang \cite{zhang2019discovering}.

\subsubsection{Strategies for method construction}
In related forthcoming work, we have discovered structure similar to the properties
discussed above among other RK methods of very high order, including some 16-stage
10th-order explicit RK methods found by numerical search \cite{zhang2019discovering}.
This suggests that these kinds of properties may be necessary for the construction of
high-order methods with low stage count, or may at least be typical of most methods in such classes.

This can be used to simplify the search for such methods by either numerical or analytical means.
One can begin by assuming the stage orders of the stages (with most stages having a selected higher order
and one or a few having lower order), and constraining the corresponding weights $b_j$ to vanish.
One can then neglect a potentially large fraction of the order conditions (as they will be
automatically satisfied) and proceed with numerical or analytical solution of those that remain.
Since it is not known \emph{a priori} which stages should have lower or higher order, the problem
is converted into a mixed-integer program with a very small number of integer variables but a
drastically reduced set of constraints and a slightly reduced set of continuous variables.

The power of this strategy remains to be explored and is beyond the scope of the present work, but
we remark that it becomes increasingly significant as the design order increases.  The structural
properties of the methods from \cite{zhang2019discovering} make it possible to reduce the
1205 order conditions to less than 400, after imposing just a few stage-order conditions and forcing
certain weights to vanish.

\section{Numerical Tests}\label{sec5}
In the following, we test the performance of the new DIRK-type schemes on several stiff, singularly perturbed, and differential-algebraic systems. When available, the analytic solution was used to compute the approximation errors. Otherwise, a highly accurate numerical reference solution was used.

\subsection{Linear Stiff ODE: Prothero-Robinson Problem}
To characterize the stability and accuracy of numerical methods for stiff problems, Prothero and Robinson \cite{prothero1974stability} introduced a family of ODEs of the form
\begin{align} \label{PR}
    y'(t) = \mu (y(t) - g(t)) + g'(t),  \ \ \ \ y(0) = g(0), \ \ \ \ \RE(\mu) \leq 0.
\end{align}
where $g(t)$ is some smooth bounded function with moderately varying rate and $\mu$ is a stiffness parameter. For any $\mu \in \CC$, the analytic solution for this problem with the chosen initial condition is $y(t) = g(t)$. The problem is stiff if $\RE(\mu) \ll -1$. It was shown in \cite{prothero1974stability} that the application of an A-stable method to $\eqref{PR}$ may result in a larger error in the stiff regime (i.e., when $\Delta t \cdot \RE(-\mu) \to \infty$) in comparison with stiffly accurate methods. Moreover, it was shown that the application of implicit one-step methods (stiffly accurate or not) to integrate $\eqref{PR}$ results in a reduction in the rate of convergence of the error, predominantly for stepsizes larger than $1/\RE(-\mu)$, which represents the fastest timescale of the dynamics in the problem.  

\begin{figure}[!ht!]
\includegraphics[width=6in]{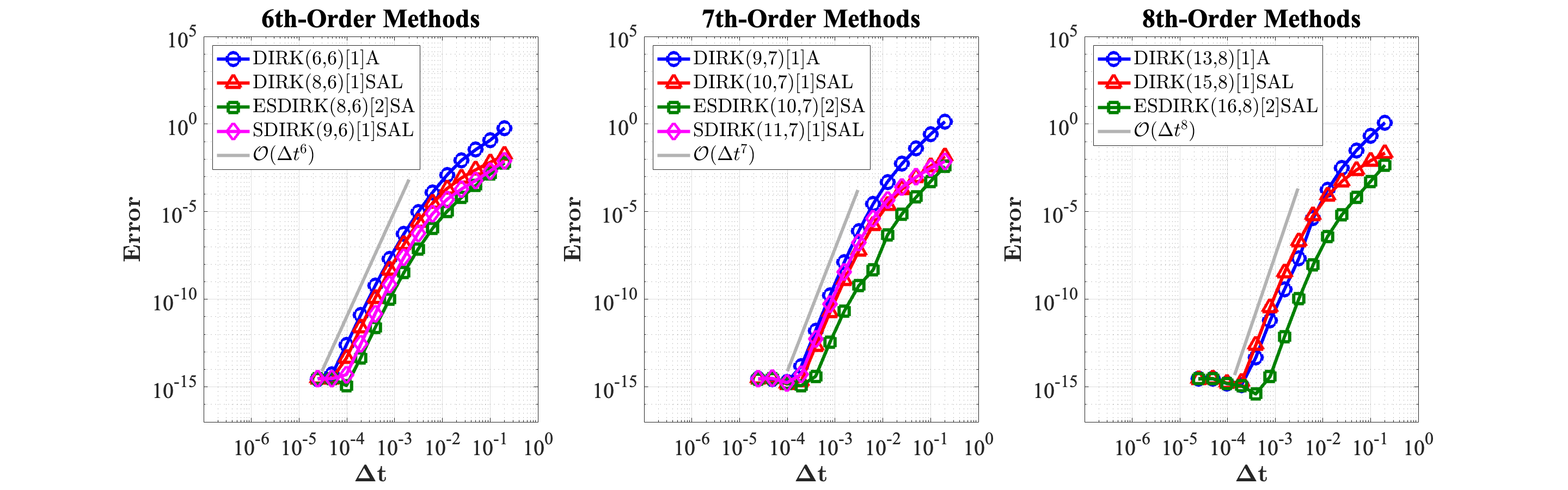}
\caption{The convergence curves for each of the 6th- (left), 7th- (middle), and 8th- (right) order schemes used to integrate the model problem $\eqref{PR}$ with stiffness parameter $\mu = -1000$.}
\label{fig:PRplt}
\end{figure}

The new methods were used to integrate $\eqref{PR}$ with a range of stepsizes $\Delta t$ over the time interval $ t \in [0,1]$ with $g(t) = \exp(-t) \cos(20t) + \sin(10t)$ and the stiffness parameter was set as $\mu = -1000$. The error was computed as the difference between the approximate and analytic solution in the maximum norm. Figure~\ref{fig:PRplt} presents the convergence curves for each of the 6th- (left), 7th- (middle), and 8th- (right) order schemes. Reference curves of the corresponding slopes are added for comparison in gray. For stepsizes smaller than $10^{-4}$, the error is around $10^{-15}$ and at that point is dominated by rounding errors. For $\Delta t \in [10^{-4},10^{-3}]$, the convergence rate is close to the design order of the method used to integrate the ODE. A reduction in the convergence rate is observed starting at $\Delta t = 10^{-3} = 1/\RE(-\mu)$ and the severity of the reduction increases as $\Delta t$ increases. Notably, the stage order 2 methods (green curves) result in the lowest error magnitude per stepsize, and along with the solely A-stable methods (blue curves) show least reduction in convergence rate in the stiff regime.

 \subsection{Singularly Perturbed Equations: Van der Pol System}
 A common class of problems for testing stiff integrators is singularly perturbed systems of the form \eqref{SPP_general}. This class of problems may help detecting the shortcomings of a given RK method, and \textit{"should be considered as a worst-case situation when the methods are applied to practical problems"}\cite{kennedy2019higherAdd}. Following the convergence study of the sixth-order nine-stage ESDIRK method presented in \cite[Fig.10]{kennedy2019diagonally}, we consider the well-known Van der Pol (VdP) system:
  \begin{subequations} \label{SPP}
\begin{align}
        y'(t) &= z, \ \ \ \ \ \ \ \ \ \ \ \ \ \ \ \ \ \ \ \  y(0) = 2 \\
        \varepsilon z'(t) &= (1-y^2)z-y, \ \ \ \ z(0) = -\dfrac{2}{3} + \dfrac{10}{81}\varepsilon - \dfrac{292}{2187}\varepsilon^2 + \dfrac{15266}{59049}\varepsilon^3 + \mathcal{O}(\varepsilon^4)
\end{align}
 \end{subequations}
 for $0 \leq t \leq 0.5$ (see also \cite{gonzalez2011global} and \cite[Sec. VI.3]{wanner1996solving}). No exact solution for the VdP system is available; hence, the reference solution we employ is obtained by integrating with the MATLAB's \texttt{ode89} solver (based on an 8th order explicit RK method) where the maximum stepsize was set to $\Delta t = 10^{-7}$ and the relative and absolute tolerances were set to $10^{-14}$ and $10^{-16}$, respectively. The stiffness of this system increases as $\varepsilon \to 0$, and if $\varepsilon = 0$, it becomes a semi-explicit index-1 DAE.  
 
 When applying an A-stable stiffly accurate $p$th order, stage order $r$ RK method to $\eqref{SPP}$, the expected convergence rate for $y$ is $\mathcal{O}\left((\Delta t)^p \right) + \mathcal{O}\left(\varepsilon (\Delta t)^{r+1}\right)$ and that of $z$ is $\mathcal{O}\left((\Delta t)^p \right) + \mathcal{O}\left(\varepsilon (\Delta t)^{r}\right)$ \cite[Sec. VI. 3, Corollary 3.10]{wanner1996solving}. Hence, in the stiff regime (when $\varepsilon \ll \Delta t$), the convergence rates for the two components decrease from the classical order $p$ to at worst $r$ and $r+1$.

 The convergence rates for the two components of the solution to \eqref{SPP} were estimated using least-square fitting \cite{kennedy2019diagonally}. The stiffness levels tested a range from the non-stiff $\varepsilon = 1$ to the highly stiff $\varepsilon = 10^{-8}$ regimes. For each stiffness level, the maximum errors for a range of stepsizes $\Delta t \in [0.25, 10^{-5}]$ were computed. All parts of each convergence curve, except the part where it levels off as $\Delta t \to 0$,  were used to infer the convergence rate, including the parts where there is an observed order reduction. 

\begin{figure}[!ht!]
\includegraphics[width=6in]{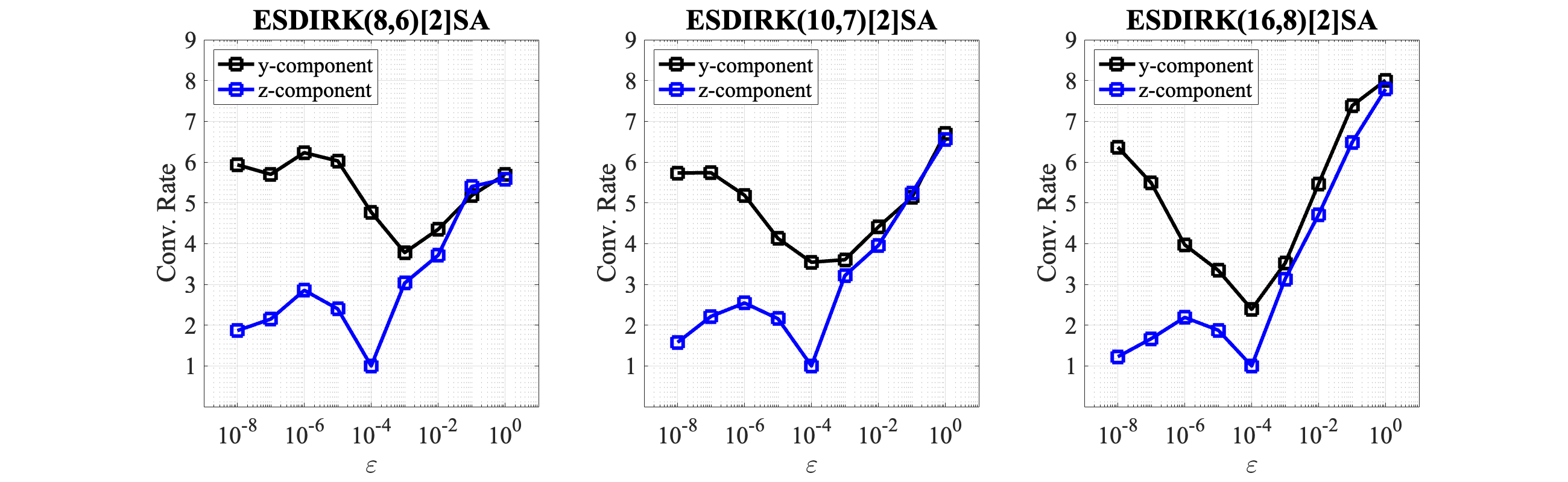}
\caption{The estimated convergence rate as a function of the stiffness parameter $\varepsilon$ for the 6th- (left), 7th- (middle), and 8th- (right) order ESDIRK schemes for the $y$-component (differential) and $z$-component (algebraic)
in the VdP system $\eqref{SPP}$.}
\label{fig:sppPlt}
\end{figure}

In Figure~\ref{fig:sppPlt}, we show the estimated convergence rates for the solution to $\eqref{SPP}$ using the $6$th, $7$th, and $8$th order ESDIRK methods. The reduction in the order of convergence is most severe for the emergent algebraic variable (i.e., the $z$-component) and attains its minimum at $\varepsilon = 10^{-4}$. While the differential variable (i.e., the $y$-component) recovers the order of convergence as $\varepsilon \to 0$, the algebraic variable shows a transient increase followed a decrease in the order. 
Interestingly, the higher order methods exhibit slightly lower convergence rates (compared to lower-order methods) at $\varepsilon = 10^{-8}$.

\subsection{Differential-Algebraic Equation: Van der Pol System}  \label{vdpDAEsec}
\begin{figure}[!ht!]
\includegraphics[width=6.0in]{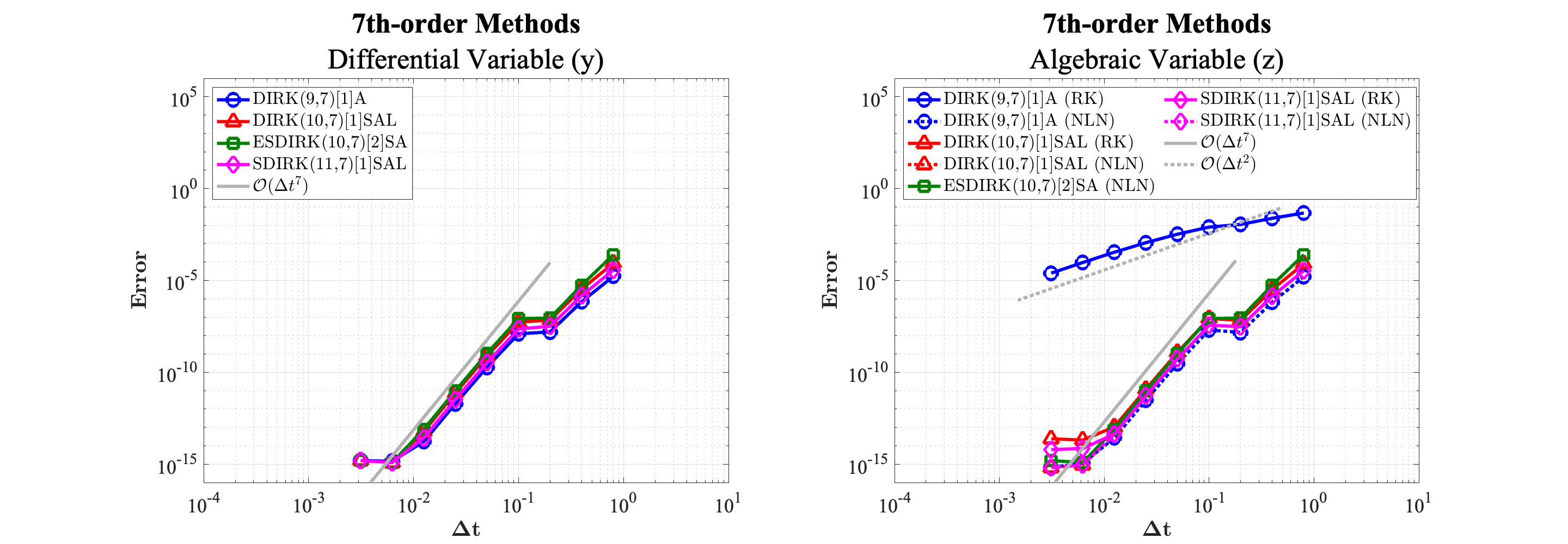}
\caption{The convergence curves for the differential variable $y$ (left) and algebraic variable $z$ (right) for the 7th order DIRK-type schemes. The algebraic variable was computed via the RK iteration \eqref{algebRK1} (RK) and by solving the nonlinear constraint \eqref{algebRK1} (NLN) \eqref{vdpDAE}.}
\label{fig:vdpDAEplt}
\end{figure}

A classical index-1 DAE is obtained from the VdP system $\eqref{SPP}$ after re-writting it in Li\'enard's coordinates and taking $\varepsilon=0$ to obtain \cite[Sec. VI.]{wanner1996solving}
\begin{subequations}\label{vdpDAE}
\begin{align}
    y' &= -z, \ \ \ y(0) = 2 \label{vdpDAE_1} \\
    0  &= y - \left( \dfrac{z^3}{3} - z \right). \label{vdpDAE_2}
\end{align}
\end{subequations}
The solution to this DAE can be approximated by an RK method with the iteration \eqref{RK_DAE}. The analytic solution for the algebraic variable is $\ln{|z|}-z^2/2 = t + K$, where $K$ is a constant that depends on the initial conditions, and the analytic solution for the differential variable can be computed using the algebraic constraint $\eqref{vdpDAE_2}$. When applying a $p$th order A-stable stiffly accurate RK method, the rate of convergence for the algebraic $z$ and differential $y$ variables are expected to be $\mathcal{O}(\Delta t^p)$. If a merely A- or L-stable RK method is used, the rate of convergence for the algebraic variable is expected to be less than $p$ and is also dependent on the stage order of the method \cite[Sec. VI.1, Theorem 1.1]{wanner1996solving}. This is demonstrated in Figure~\ref{fig:vdpDAEplt} where the convergence curves for the differential and algebraic variables are computed separately using the 7th-order DIRK-type methods as a sample. The rest of the methods exhibit analogous results.  The error was computed as the maximum norm of the difference between the analytic and approximate solution. The (consistent) chosen initial conditions are $(y(0),z(0)) = (1,2.1)$ and the DAE was integrated over the interval $t \in [0,0.9]$. The algebraic variable $z$ was approximated using two approaches: one by solving the nonlinear equation $\eqref{vdpDAE_2}$ (denoted by NLN) at each step after computing the differential variable and the other by using the RK iteration \eqref{algebRK1} (denoted by RK). Notice that the second approach for computing the algebraic variable is not possible for ESDIRK methods as $\AAA$ is singular.

As shown in the left panel of Figure~\ref{fig:vdpDAEplt}, with all of the applied 7th-order methods, the rate of convergence for the differential variable is consistent with order of the method. On the right panel, computing the algebraic variable using the nonlinear relation $\eqref{algebRK2}$ (dashed curves) results in a 7th order convergence since a sufficiently accurate nonlinear solver is used to compute $z_{n+1}$ and $y_{n+1}$ is computed with a 7th order of accuracy. Moreover, computing the algebraic variable using the RK iteration (solid curves) with the stiffly accurate methods also maintains the 7th order of convergence. Consistent with \cite[Sec. VI.1, Theorem 1.1]{wanner1996solving}, computing the algebraic variable with a stage order 1 A-stable method is expected to result in a 2nd order of convergence, which is the case with the $\texttt{DIRK(9,7)[1]A}$ method (blue solid curve). 

\subsection{Adaptive Stepsize: PID Controller} 
\begin{figure}[h]
\includegraphics[width=6in]{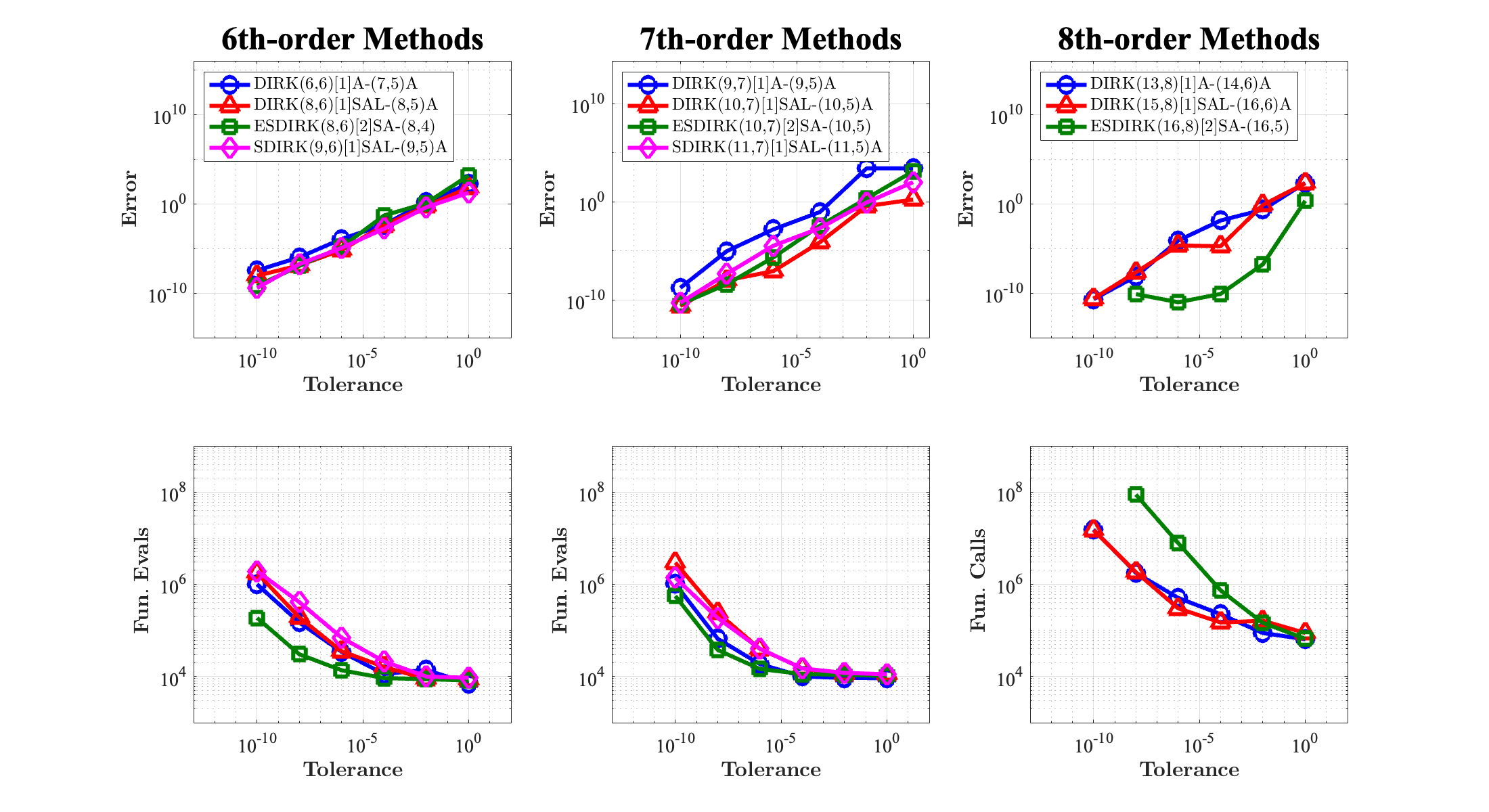}
\caption{Error vs. tolerance (top row) and the number of function calls (i.e., Newton's iterations) vs. tolerance (bottom row) for the  6th- (left), 7th- (middle), and 8th- (right) order embedded schemes applied to the VdP system $\eqref{SPP}$ with $\varepsilon = 10^{-5}$ over $t \in [0,2]$ and with PID-controller.}
\label{fig:wpPlt}
\end{figure}

Here, we demonstrate the performance of the new schemes when using the embedded error estimators along with a stepsize control routine. We approximate the solution to the VdP system $\eqref{SPP}$ with $\varepsilon = 10^{-5}$ over the interval $t = [0, 2]$. The controller used here is the Proportional-Integral-Derivative (PID) controller (see e.g., \cite{soderlind2003digital}, \cite{wanner1996solving}]) where the new time step is updated as:
$$\Delta t_{n+1} = w_{n+1}^{-\beta_1/(\widehat{p}+1)} w_{n}^{-\beta_2/(\widehat{p}+1)} \Delta t_n$$
with
$$w_n := \sqrt{ \dfrac{1}{m} \sum_{i=1}^{m} \left( \dfrac{y^i_{n} - \hat{y}^i_n}{\texttt{rtol}^i\cdot \max\{|y^i_n|,|\hat{y}^{i}_n|\} + \texttt{atol}^i} \right)^2 }$$
$\widehat{p}$  is the order of the error estimator, $w_{n}, w_{n+1}$ are scaled local errors, and $\texttt{atol}, \texttt{rtol}$ are absolute and relative error tolerances (here we take $\texttt{atol} = \texttt{rtol}$ for all solution components). The subscripts represent the index of the time step and the superscripts denote the component of the solutions (for the VdP, $m=2$). The exponents are chosen as $\beta_1 = 0.60$ and $\beta_2=-0.20$ (in the related literature, this is usually identified as PI42 controller). The initial stepsize is  $\Delta t_0 = 10^{-8}$. The measures of efficiency we have chosen to present in Figure.~\ref{fig:wpPlt} are the global error, computed as the root mean-squared (RMS) error, and the number of function calls (i.e., number of Newton's iterations) for the right hand side vs. the imposed relative (and accordingly the absolute) tolerances.

It can be observed that with all of the embedded pairs (with the exception of the \texttt{ESDIRK(16,8)[2]SAL-[(16,5)]} scheme where the incurred RMS error is lower than what is requested), the incurred RMS error is approximately the imposed local error tolerance.  Whereas one would typically expect higher-order methods to be more efficient, this is not necessarily the case here due to order reduction.

\subsection{Stiff PDE Semidiscretization: The Heat Equation}
\begin{figure}[!ht!]
\includegraphics[width=6in]{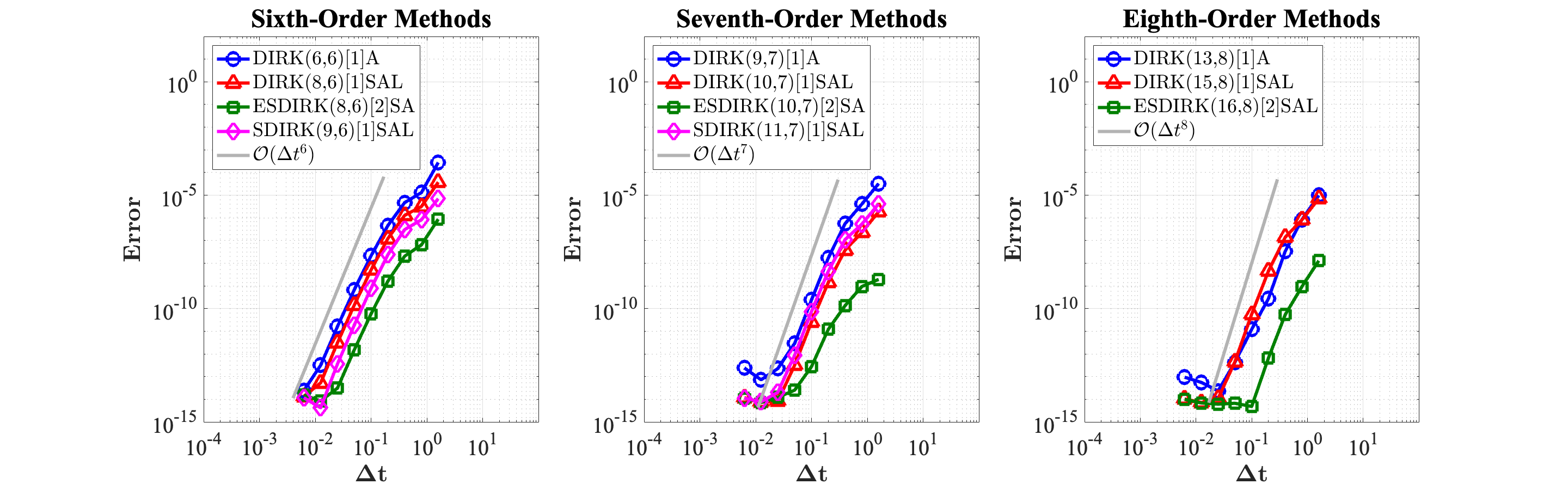}
\caption{The convergence curves for each of the 6th- (left), 7th- (middle), and 8th- (right) order schemes used to integrate the semidiscretization of $\eqref{heat}$.}
\label{fig:HEAT}
\end{figure}

Larger stiff systems of ODEs may arise in method of lines discretizations of PDEs. As an illustration, we consider the numerical approximation of the solution to the one-dimensional heat equation 
\begin{equation}\label{heat}
    u_t - u_{xx} = g(x,t); \ \ \ \ (x,t) \in (0,1) \times [0,5]
\end{equation}
where $g(x,t)$ and the initial and boundary conditions are chosen such that $u(x,t) = \exp(-t/10)\sin(\pi x)$. The second-order centered finite difference scheme was used to discretize $u_{xx}$ in space which resulted in an $m$-dimensional semidiscrete system of ODEs whose eigenvalues reside on the negative real axis.  The stiffness of such a system (which is determined by the magnitude of the eigenvalues) increases as the spatial grid size approaches 0. The semidiscrete system was integrated with the presented schemes with $m=200$ and $\Delta x = 1/m$.  The error was computed as the maximum norm in space of the difference between the approximate solution and the solution obtained with MATLAB's stiff solver $\texttt{ode15s}$ at $t = 5$ (with the relative tolerance \texttt{RelTol} in was set to $10^{-14}$ and absolute tolerance \texttt{AbsTol} was set $10^{-16}$). The convergence curves for the 6th, 7th, and 8th order schemes are presented in figure~\ref{fig:HEAT}. 
As one might expect, the best accuracy is obtained with methods of stage order two.  The methods achieve their design order of convergence over a range of step sizes, although this range is smaller for very high order methods since they reach the level of roundoff error very quickly.

\section{Conclusions}\label{sec6}
In this work, we have presented new sixth, seventh, and eighth order A-stable DIRK-type schemes, some of which are additionally L-stable, stiffly accurate, and/or have stage order two. We found these schemes by numerical search based on combined global and local optimization strategies where the constraints were formulated directly in terms of the coefficients of the schemes and certain reformulations of the original optimization problem were employed. We have described several approaches that enabled finding highly accurate numerical solutions to the constraints introduced by the high-order polynomials associated with the RK order conditions and the A-stability requirement. The approach presented herein can be readily adopted to search for various types of time integrators with different stability and structure properties. Such searches may be aided by incorporating some of the
simplifying structural properties discussed in Section~\ref{SP}.

The performance of the schemes was demonstrated on a variety of (moderately and highly) stiff and differential-algebraic systems of ODEs and semidiscretizations of PDEs. As expected when integrating with implicit schemes, the order reduction phenomenon was observed and its severity is most pronounced in the stiff regime of a given problem (Figure~\ref{fig:PRplt}) or when the controlling parameter is chosen such that the problem is highly stiff (Figure~\ref{fig:sppPlt}). 
Besides, for highly stiff problem, the severity of the order reduction increases with the order of the applied method. Taking into account these observations, the choice of an appropriate DIRK-type method is problem-dependent.

However, as a general guideline, since the presented DIRK-type methods are designed with no targeted class of problems, they may serve as a good choice for integrating modestly stiff systems where the constraint on the stepsize required by explicit methods needs to be avoided. Employing the embedded error estimators combined with suitable stepsize control routines can further reduce the computational cost incurred when using a constant stepsize value for entire domain of integration. Among all of the presented methods, higher order methods can be used to achieve lower error magnitudes at a given stepsize value. Except in the approximation of algebraic variables in DAEs (e.g., Figure~\ref{fig:vdpDAEplt}), stage order one, L-stable, and stiffly accurate methods do not demonstrate a noticeable advantage over their solely A-stable counterparts (except yielding slightly lower error magnitudes in some problems). The former group of methods will be advantageous if damping of fast transients in the numerical solution is needed. On the other hand, stage order two methods, despite having more stages, result in higher accuracy overall in comparison with the rest of the methods. Since the stage order for DIRK-type schemes is restricted to 2 or less, a potential promising future direction is the design of similar schemes with high \emph{weak stage order} \cite{2018_wso}, an approach that have shown promising results in avoiding order reduction for certain classes of stiff problems.

Several other future directions based on the techniques developed in this work are possible. One is to further enhance the presented methods with additional capabilities that improve their efficiency and practicality such as equipping them with dense-output formulas and optimizing them to facilitate the use of sophisticated stage-value predictors. Another is the design of high-order methods optimized for specific classes of problems or applications, such as designing methods with optimized dispersion and dissipation errors. Finally, global search strategies may be employed to aid in establishing the true minimum number of stages for high order explicit and DIRK methods.

\section*{Acknowledgement} The computations presented in this work were performed in the IBEX clusters at King
Abdullah University of Science and Technology (KAUST). The authors wish to thank Ricardo de Lima for various insightful discussions that enabled solving the optimization problems presented in this work, Steven Roberts for suggesting the use of local optimization in Mathematica to improve the accuracy of the presented schemes, and an anonymous referee whose suggestions led to useful extensions and improvements of this work.

\newpage
\section*{Appendix: Coefficients of the DIRK-type schemes}

\begin{enumerate} [1.]
    \item The \texttt{DIRK(6,6)[1]A-[(7,5)A]} scheme.
    \centering
     \begingroup
\ttfamily

\begin{tabular}{%
>{\begin{math}}c<{\end{math}}r 
}
\midrule
\multicolumn{1}{c}{\textrm{\textbf{Coefficient}}} & \multicolumn{1}{c}{\textrm{\textbf{Value}}} \\
\toprule
a_{1,1}  &  3.034878447067473e-01   \\
a_{2,1}  &   -2.797564927098137e-01  \\
a_{2,2}  &  5.000322360207473e-01  \\
a_{3,1}  &  2.805832157438948e-01         \\
a_{3,2}  & -4.385600615867512e-01    \\
a_{3,3}  &  2.172507345157362e-01  \\
a_{4,1}  &   -6.776787385398458e-02              \\
\bottomrule
\end{tabular}
\endgroup

 \newpage
     \item The \texttt{DIRK(8,6)[1]SAL-[(8,5)A]} scheme.
         \centering
     \begingroup
\ttfamily

\begin{tabular}{%
>{\begin{math}}c<{\end{math}}r 
}
\midrule
\multicolumn{1}{c}{\textrm{\textbf{Coefficient}}} & \multicolumn{1}{c}{\textrm{\textbf{Value}}} \\
\toprule
a_{1,1}  &  4.772644573858262e-01   \\
a_{2,1}  &   -1.970525884150017e-01     \\
a_{2,2}  &   4.763634284595835e-01     \\
a_{3,1}  &  -3.476744303729656e-02       \\
a_{3,2}  &   6.330518073354831e-01     \\
a_{3,3}  &    1.936343100750279e-01     \\
a_{4,1}  &  9.677976685787021e-02          \\
\bottomrule
\end{tabular}
\endgroup

    \newpage
     \item The \texttt{ESDIRK(8,6)[2]SA-[(8,4)]} scheme.
         \centering
     \begingroup
\ttfamily

\begin{tabular}{%
>{\begin{math}}c<{\end{math}}r 
}
\midrule
\multicolumn{1}{c}{\textrm{\textbf{Coefficient}}} & \multicolumn{1}{c}{\textrm{\textbf{Value}}} \\
\toprule
a_{1,1}  &  0.000000000000000e+00   \\
a_{2,1}  = a_{i,i}  = b_8 &   3.332221492177252e-01     \\
a_{3,1}  &    6.397437731822139e-02       \\
a_{3,2}  &    -8.303302244102144e-02    \\
a_{4,1}  &     -7.285222013693263e-01      \\
a_{4,2} &   -2.104144795224848e-01   \\
a_{4,3} &     5.325199165593416e-01    \\
\bottomrule
\end{tabular}
\endgroup

    \newpage
     \item The \texttt{SDIRK(9,6)[1]SAL-[(9,5)A]} scheme.
         \centering
     \begingroup
\ttfamily

\begin{tabular}{%
>{\begin{math}}c<{\end{math}}r 
}
\midrule
\multicolumn{1}{c}{\textrm{\textbf{Coefficient}}} & \multicolumn{1}{c}{\textrm{\textbf{Value}}} \\
\toprule
a_{i,i} = b_9 &  2.181277819449076e-01   \\
a_{2,1}  &    -9.035148561194185e-02        \\
a_{3,1}  &    1.729520391389366e-01    \\
a_{3,2}  &   -3.536550103628203e-01     \\
a_{4,1}  &       5.119998759191926e-01    \\
a_{4,2} &     2.896403322019248e-02    \\
a_{4,3} &   -1.440309456570937e-02    \\
\bottomrule
\end{tabular}
\endgroup

    \newpage
    \item The \texttt{DIRK(9,7)[1]A-[(9,5)A]} scheme.
    \centering
     \begingroup
\ttfamily
\begin{tabular}{%
>{\begin{math}}c<{\end{math}}r 
>{\begin{math}}c<{\end{math}}r 
}
\midrule
\multicolumn{1}{c}{\textrm{\textbf{Coefficient}}} & \multicolumn{1}{c}{\textrm{\textbf{Value}}} & \multicolumn{1}{c}{\textrm{\textbf{Coefficient}}} & \multicolumn{1}{c}{\textrm{\textbf{Value}}} \\
\toprule
a_{1,1} &  1.798777898558387e-01    &  a_{8,4} &  -1.348926341392102e-01                      \\
a_{2,1}  &    -1.004058448851567e-01       &  a_{8,5} & 1.295532472606773e-01   \\ 
a_{2,2}  &     2.149485906448186e-01       &  a_{8,6} &   9.924937953714891e-02     \\
a_{3,1}  &   1.122513601989947e-01         &  a_{8,7} &   -1.571661056346096e-01   \\
a_{3,2} &   -2.061621391502978e-01      &  a_{8,8} &    1.791886227981400e-01   \\
a_{3,3} &     1.251596429419579e-01    &   a_{9,1} &  4.396960799652251e-03    \\
a_{4,1}  &     -3.351640007682573e-02     &  a_{9,2} &   9.602504865704912e-01     \\
\bottomrule
\end{tabular}
\endgroup

     \newpage
    \item The \texttt{DIRK(10,7)[1]SAL-[(10,5)A]} scheme.
    
\centering
     \begingroup
\ttfamily
\begin{tabular}{%
>{\begin{math}}c<{\end{math}}r 
>{\begin{math}}c<{\end{math}}r 
}
\midrule
\multicolumn{1}{c}{\textrm{\textbf{Coefficient}}} & \multicolumn{1}{c}{\textrm{\textbf{Value}}} & \multicolumn{1}{c}{\textrm{\textbf{Coefficient}}} & \multicolumn{1}{c}{\textrm{\textbf{Value}}} \\
\toprule
a_{1,1} &  2.337046321252643e-01     &  a_{8,5} &   -1.617554588390476e-04   \\
a_{2,1}  &  -7.393248131494073e-02    &  a_{8,6} &   -8.481343359802809e-02  \\ 
a_{2,2}  &   2.000568381461040e-01   &  a_{8,7} &    -3.643866656666597e-01   \\
a_{3,1}  &  9.437903440448119e-02         &  a_{8,8} &  1.586044201360552e-01   \\
a_{3,2} &    2.640560677016047e-01     &  a_{9,1} &        2.007331564774248e-01 \\
a_{3,3} &   1.332452024564650e-01   &   a_{9,2} &     2.396864434444326e-01   \\
a_{4,1}  &  2.690848106012005e-01           &  a_{9,3} & 3.038370144189287e-01  \\
\bottomrule
\end{tabular}
\endgroup

     \newpage
    \item The \texttt{ESDIRK(10,7)[2]SA-[(10,5)]} scheme.
    
\centering
     \begingroup
\ttfamily
\begin{tabular}{%
>{\begin{math}}c<{\end{math}}r 
>{\begin{math}}c<{\end{math}}r 
}
\midrule
\multicolumn{1}{c}{\textrm{\textbf{Coefficient}}} & \multicolumn{1}{c}{\textrm{\textbf{Value}}} & \multicolumn{1}{c}{\textrm{\textbf{Coefficient}}} & \multicolumn{1}{c}{\textrm{\textbf{Value}}} \\
\toprule
a_{1,1} &  0.000000000000000e+00     &  a_{8,7} &     1.592795834073081e-02       \\
a_{2,1} = a_{i,i} = b_{10} &  2.100557902034192e-01   &  a_{9,1} &   -3.030621290769446e-02    \\ 
a_{3,1}  &   2.557817399210862e-01   &  a_{9,2} &   -2.970351746590337e-01    \\
a_{3,2}  &     2.398509169809763e-01       &  a_{9,3} &   1.847246974621637e-01    \\
a_{4,1} &   2.867896248804369e-01    &  a_{9,4} &    -3.518760795161831e-02      \\
a_{4,2} &    2.304947488347783e-01   &   a_{9,5} &   -3.246682306907609e-03     \\
a_{4,3} &    2.639251498854907e-01        &  a_{9,6} &   2.161510040535314e-01    \\
\bottomrule
\end{tabular}
\endgroup

     \newpage
    \item The \texttt{SDIRK(11,7)[1]SAL-[(11,5)A]} scheme.
    
\centering
     \begingroup
\ttfamily
\begin{tabular}{%
>{\begin{math}}c<{\end{math}}r 
>{\begin{math}}c<{\end{math}}r 
}
\midrule
\multicolumn{1}{c}{\textrm{\textbf{Coefficient}}} & \multicolumn{1}{c}{\textrm{\textbf{Value}}} & \multicolumn{1}{c}{\textrm{\textbf{Coefficient}}} & \multicolumn{1}{c}{\textrm{\textbf{Value}}} \\
\toprule
a_{i,i} =   b_{11} &   2.002526611877419e-01      &  a_{9,5} &    -1.176286417178894e-02     \\
a_{2,1} &  -8.294736816526700e-02   &  a_{9,6} &   4.321768808679651e-02     \\ 
a_{3,1} &   4.834526905407508e-01    &  a_{9,7} &   -3.152068362754727e-02    \\
a_{3,2} &    0.000000000000000e+00         &  a_{9,8} &   -8.460070216387972e-02     \\
a_{4,1} &   7.710764534813213e-01       &  a_{10,1} &       6.514285986237711e-01  \\                      
a_{4,2} &  -2.293692634184203e-01    &   a_{10,2} &   -1.020807847535596e-01      \\
a_{4,3} &   2.897333732088233e-01      &  a_{10,3} &    1.983057018018881e-01    \\
\bottomrule
\end{tabular}
\endgroup

    \newpage
    \item The \texttt{DIRK(13,8)[1]A-[(14,6)A]} scheme.
 
 \centering
     \begingroup
\ttfamily
\begin{tabular}{%
>{\begin{math}}c<{\end{math}}r 
>{\begin{math}}c<{\end{math}}r 
}
\midrule
\multicolumn{1}{c}{\textrm{\textbf{Coefficient}}} & \multicolumn{1}{c}{\textrm{\textbf{Value}}} & \multicolumn{1}{c}{\textrm{\textbf{Coefficient}}} & \multicolumn{1}{c}{\textrm{\textbf{Value}}} \\
\toprule
a_{1,1} &  4.210507454422905e-01    &  a_{11,5} &   -1.545062646028905e-02   \\
a_{2,1}  & -7.610794195912682e-02     &  a_{11,6} & -4.622465919227501e-01   \\ 
a_{2,2}  &   2.643539865808565e-01      &  a_{11,7} & -5.764063273291814e-02   \\
a_{3,1}  &  7.271069041706940e-02      &  a_{11,8} & -7.120669425040176e-03  \\
a_{3,2} &  -2.042659769772850e-01      &  a_{11,9} &  3.777765580144523e-01   \\
a_{3,3} &    1.816081965441359e-01     &   a_{11,10} &  3.689005433829403e-01    \\
a_{4,1}  &   5.576305481661099e-01    &  a_{11,11} & 6.184887463318368e-02  \\
\bottomrule
\end{tabular}
\endgroup

\newpage
    
 \centering
     \begingroup
\ttfamily
\begin{tabular}{%
>{\begin{math}}c<{\end{math}}r 
>{\begin{math}}c<{\end{math}}r 
}
\midrule
\multicolumn{1}{c}{\textrm{\textbf{Coefficient}}} & \multicolumn{1}{c}{\textrm{\textbf{Value}}} & \multicolumn{1}{c}{\textrm{\textbf{Coefficient}}} & \multicolumn{1}{c}{\textrm{\textbf{Value}}} \\
\toprule
a_{9,4}  &   -3.174588332428040e-03    &  a_{14,8} & -1.899130199303695e-02  \\
a_{9,5}  & -1.128903821351927e-01  & a_{14,9} & 2.206631140820358e-01  \\
a_{9,6}  & -4.198092679542844e-01    & a_{14,10} &  2.530299843608644e-01  \\
a_{9,7}  &  4.661369021021040e-02     & a_{14,11} &  2.520117993705631e-01   \\
a_{9,8}  &   -1.154138130410853e-03    &  a_{14,12} &   -1.548822226054231e-01   \\
a_{9,9}  &      1.096853636923826e-01     &  a_{14,13} & 3.152022646874153e-02  \\
a_{10,1}  &     -5.126168782523551e-02  &  a_{14,14} &   5.140958121047144e-02       \\
\bottomrule
\end{tabular}
\endgroup

    \newpage
    \item The \texttt{DIRK(15,8)[1]SAL-[(16,6)A]} scheme.
    
    \centering
     \begingroup
\ttfamily
\begin{tabular}{%
>{\begin{math}}c<{\end{math}}r 
>{\begin{math}}c<{\end{math}}r 
}
\midrule
\multicolumn{1}{c}{\textrm{\textbf{Coefficient}}} & \multicolumn{1}{c}{\textrm{\textbf{Value}}} & \multicolumn{1}{c}{\textrm{\textbf{Coefficient}}} & \multicolumn{1}{c}{\textrm{\textbf{Value}}} \\
\toprule
a_{1,1} & 4.989049812711929e-01    &  a_{9,4} &  -1.822496722358612e-01    \\
a_{2,1}  &  -3.038060373418158e-01      &  a_{9,5} &  1.678976357137985e-01   \\ 
a_{2,2}  &  8.862994459923792e-01    &  a_{9,6} & 2.128503350300694e-01    \\
a_{3,1}  &  -5.814402234714763e-01          &  a_{9,7} &  -3.917392994401228e-01  \\
a_{3,2} &   3.710037194602585e-01       &  a_{9,8} &   -1.187185068767668e-02    \\
a_{3,3} &     4.384471775280196e-01    &   a_{9,9} &   5.262937016590934e-01   \\
a_{4,1}  &   5.318526388700512e-01            &  a_{10,1} &      3.839839148454607e-01        \\
\bottomrule
\end{tabular}
\endgroup

\newpage
    
 \centering
     \begingroup
\ttfamily
\begin{tabular}{%
>{\begin{math}}c<{\end{math}}r 
>{\begin{math}}c<{\end{math}}r 
}
\midrule
\multicolumn{1}{c}{\textrm{\textbf{Coefficient}}} & \multicolumn{1}{c}{\textrm{\textbf{Value}}} & \multicolumn{1}{c}{\textrm{\textbf{Coefficient}}} & \multicolumn{1}{c}{\textrm{\textbf{Value}}} \\
\toprule
a_{13,1}  &      1.570142745612990e-01  & a_{15,10} = b_{10}    &    7.982503922188518e-02  \\
a_{13,2} &   -8.838102563818741e-02   & a_{15,11}  = b_{11}    &     5.288247340826546e-02  \\
a_{13,3}  & 1.171930338850341e-01    & a_{15,12} = b_{12}  &    -8.303508889003619e-02   \\
a_{13,4}  &     -3.623042437694662e-02    & a_{15,13} = b_{13}   &    2.256738870727901e-02 \\
a_{13,5}  & 1.690302114661112e-02     & a_{15,14} = b_{14}  &   -5.926311190402041e-02  \\
a_{13,6}  &  -1.698357535761414e-01     & a_{15,15} = b_{15}  &  1.068258780376210e-01   \\
a_{13,7}  &       3.997499792341125e-01   & a_{16,1}     &     1.737844812076520e-01    \\
\bottomrule
\end{tabular}
\endgroup

    \newpage
    \item The \texttt{ESDIRK(16,8)[2]SAL-[(16,5)]} scheme.
    
    \centering
     \begingroup
\ttfamily
\begin{tabular}{%
>{\begin{math}}c<{\end{math}}r 
>{\begin{math}}c<{\end{math}}r 
}
\midrule
\multicolumn{1}{c}{\textrm{\textbf{Coefficient}}} & \multicolumn{1}{c}{\textrm{\textbf{Value}}} & \multicolumn{1}{c}{\textrm{\textbf{Coefficient}}} & \multicolumn{1}{c}{\textrm{\textbf{Value}}} \\
\toprule
a_{1,1} &  0.000000000000000e+00   &  a_{10,3} &    1.754111303661099e-02          \\
a_{2,1} = a_{i,i} = b_{16} &   1.173188193585213e-01   &  a_{10,4} &   1.448715276824183e-01     \\ 
a_{3,1}  &   5.570146059746160e-02   &  a_{10,5} &   -4.185545976980595e-02  \\
a_{3,2}  &   3.855256466387416e-01  &  a_{10,6} &       8.418321683322609e-02    \\
a_{4,1} &  6.349327642889503e-02     &  a_{10,7} &  -8.500209372821918e-02     \\
a_{4,2} &   3.735561262636805e-01    &   a_{10,8} &     4.861703438258985e-01   \\
a_{4,3}  &   8.299416643895295e-03   &  a_{10,9} &   -5.267171168227393e-02  \\
\bottomrule
\end{tabular}
\endgroup

\newpage
    
 \centering
     \begingroup
\ttfamily
\begin{tabular}{%
>{\begin{math}}c<{\end{math}}r 
>{\begin{math}}c<{\end{math}}r 
}
\midrule
\multicolumn{1}{c}{\textrm{\textbf{Coefficient}}} & \multicolumn{1}{c}{\textrm{\textbf{Value}}} & \multicolumn{1}{c}{\textrm{\textbf{Coefficient}}} & \multicolumn{1}{c}{\textrm{\textbf{Value}}} \\
\toprule
a_{13,12}  &  8.841296671144815e-02     & a_{16,2} = b_{2}    &  7.800753944828058e-02 \\
a_{14,1} &   -5.455619138481062e-02   & a_{16,3}  = b_{3}    &  1.502187414805799e-02  \\
a_{14,2}  &  6.757854234427529e-01    & a_{16,4} = b_{4}  &   1.951802772841952e-01   \\
a_{14,3}  &    4.230664432019405e-01   & a_{16,5} = b_{5}   &  -2.466433101532346e-03     \\
a_{14,4}  &  -1.653001268411932e-04     & a_{16,6} = b_{6}  &  4.739771170683139e-02  \\
a_{14,5}  &  1.042529947937628e-01    & a_{16,7} = b_{7}  &   -6.827735586103630e-02         \\
a_{14,6}  &    -1.057630193030210e-01   & a_{16,8}  = b_{8}    &  1.956801912387796e-01    \\
\bottomrule
\end{tabular}
\endgroup

\end{enumerate}

\newpage

\bibliographystyle{plain} 
\bibliography{ref}

\section*{Declarations}
\subsection*{Funding}
This work was supported by funding from King Abdullah University of Science and Technology.

\subsection*{Competing interests}
The authors have no relevant financial or non-financial interests to disclose.

\subsection*{Data Availability}
The coefficients of the methods developed in this study are available in the appendix and also online from \url{https://github.com/yousefalamri55/High_Order_DIRK_Methods_Coeffs}.
\end{document}